\documentclass{article}
\usepackage{amssymb,amsmath,theorem,a4wide}
\usepackage{graphicx}

\allowdisplaybreaks[2]

\newcommand\eq{\leftrightarrow}
\newcommand\LOR{\bigvee}
\newcommand\ET{\bigwedge}
\newcommand\model{\vDash}
\newcommand\nmodel{\nvDash}
\newcommand\eqs{\bumpeq}

\newcommand\fii{\varphi}
\newcommand\tet{\vartheta}
\newcommand\ep{\varepsilon}

\newcommand\p[1]{\langle#1\rangle}
\newcommand\lh[1]{\lvert#1\rvert}
\newcommand\vc[1]{\overline{#1}}
\newcommand\cat{^\smallfrown}
\newcommand\bez{\smallsetminus}
\newcommand\sset{\subseteq}
\newcommand\nsset{\nsubseteq}
\newcommand\ssset{\subsetneq}
\newcommand\Sset{\supseteq}
\newcommand\cupd{\mathbin{\dot\cup}}
\newcommand\nul{\varnothing}
\newcommand\res{\mathbin\restriction}
\newcommand\id{\mathrm{id}}
\newcommand\pto{\rightharpoonup}

\newcommand\code[1]{\ulcorner#1\urcorner}
\newcommand\num[1]{\underline{#1}}
\newcommand\jvn[1]{\boldsymbol{#1}}
\newcommand\sucf{\mathrm{succ}}

\DeclareMathOperator\Lh{lh}
\DeclareMathOperator\acl{acl}
\DeclareMathOperator\dcl{dcl}
\DeclareMathOperator\tp{tp}
\DeclareMathOperator\diag{Diag}

\newcommand\sprf{\mathrm{PRF}}
\newcommand\strf{\mathrm{TRF}}
\newcommand\sdprp{\mathrm{DPRP}}
\newcommand\srp{\mathrm{RP}}
\newcommand\prpty[1]{\mathrm{#1}}

\newcommand\ec{\mathit{EC}}
\newcommand\trep[1]{\mathit{REP}_{#1}}
\newcommand\tprf{\trep\sprf}
\newcommand\tprfu{\trep U}
\newcommand\teq{\mathrm{eq}}

\newcommand\N{\mathbb N}
\newcommand\monster{\mathbb M}

\newcommand\indep[1][]{\doindep\vert{#1}}
\newcommand\nindep[1][]{\doindep\nvert{#1}}
\newcommand\doindep[2]{\mathrel{\mathop{\vcenter{%
        \hbox{\oalign{\noalign{\kern-.3ex}\hfil$#1$\hfil\cr
              \noalign{\kern-.7ex}%
              $\smile$\cr\noalign{\kern-.3ex}}}}}\displaylimits_{#2}}\nobreak}
\newcommand\nvert{{\setbox0\hbox{$\smallsetminus$}\hbox to\wd0{\hss$\vert$\hss}\kern-\wd0\reflectbox{\box0}}}

\newcommand\ob[1]{\overline{#1}}
\newcommand\ub[1]{\underline{#1}}
\newcommand\txto{${}\to{}$}

\newcommand\bme{\hskip.75em\relax}
\newcommand\noproof{\leavevmode\unskip\bme\vadjust{}\nobreak\hfill$\qed$\par}
\newcommand\qed{\Box}
\newenvironment{Pf}
  {\par\noindent\textit{Proof:}\bme\ignorespaces}
  {\noproof\pagebreak[2]\vskip\medskipamount\ignorespacesafterend}

\theoremstyle{plain}
\newtheorem{Thm}{Theorem}[section]
\newtheorem{Prop}[Thm]{Proposition}
\newtheorem{Cor}[Thm]{Corollary}
\newtheorem{Lem}[Thm]{Lemma}
\newtheorem{Obs}[Thm]{Observation}
\newtheorem{Que}[Thm]{Question}
\newtheorem{Cl}{Claim}[Thm]

\theorembodyfont\upshape
\newtheorem{Def}[Thm]{Definition}
\newtheorem{Rem}[Thm]{Remark}
\newenvironment{Pf*}{\let\qed\qedCl\Pf}\endPf

\newcommand\allowhyphens{\nobreak\hskip0pt\relax}
\DeclareRobustCommand*\LP{\ifmmode(\else\textup(\allowhyphens\fi}
\DeclareRobustCommand*\RP{\ifmmode)\else\textup)\fi}
\newcommand\paren[1]{\LP#1\RP}

\author{Emil Je\v r\'abek\\[\medskipamount]
Institute of Mathematics of the Czech Academy of Sciences\\
\small \v Zitn\'a 25,
115\:67 Praha 1,
Czech Republic,
email: \texttt{jerabek@math.cas.cz}
}

\title{Recursive functions and existentially closed structures}

\begin{document}
\maketitle

\begin{abstract}
The purpose of this paper is to clarify the relationship between various conditions implying essential undecidability:
our main result is that there exists a theory~$T$ in which all partially recursive functions are representable, yet $T$
does not interpret Robinson's theory~$R$. To this end, we borrow tools from model theory---specifically, we investigate
model-theoretic properties of the model completion of the empty theory in a language with function symbols. We obtain a
certain characterization of $\exists\forall$ theories interpretable in existential theories in the process.


\end{abstract}

\section{Introduction}\label{sec:introduction}

First-order theories studied by logicians may be broadly divided in two classes. One class comprises theories of
``arithmetical strength'', such as various fragments and extensions of Peano arithmetic, or set theories. They are
distinguished by their great expressive power that, on the one hand, allows them to work with all kinds of objects from
mathematical practice in a suitable encoding (indeed, some of these theories are designed to serve as foundations for all
of mathematics, e.g., ZFC), and on the other hand, makes them subject to G\"odel's incompleteness theorems and related
phenomena. The other class are theories of ``tame'' structures, for example algebraically closed or real closed fields,
vector spaces, generic structures such as the random graph, etc. These theories have low expressive power (often
manifested in classification of definable sets stemming from partial quantifier elimination), and consequently their
models have a manageable structure of a geometric nature. Tame theories tend to be decidable.

The borderline between arithmetical and tame theories is not sharply demarcated, but one typical feature of
arithmetical theories is their \emph{essential undecidability}, meaning that all consistent extensions of the theory
are undecidable. This notion was isolated by Tarski, Mostowski, and Robinson~\cite{tmr}. This classic monograph also
includes convenient methods for proving essential undecidability of a theory~$T$, which can be viewed as stand-alone
properties implying essential undecidability. In order of increasing strength, these are:
\begin{itemize}
\item $T$ can represent all partially recursive functions (prf; see below for a precise definition).
\item $T$ can interpret Robinson's theory~$R$.
\item $T$ can interpret Robinson's arithmetic~$Q$, or equivalently, the adjunctive set theory.
\end{itemize}
An even stronger condition is that of being an \emph{sequential theory} \cite{pud:cuts,vis:seq}.

Recall that Robinson's $R$ is,
essentially, a theory axiomatizing the true $\Sigma_1$~sentences of the standard model of arithmetic~$\N$; while it
is in some ways less convenient to work with than the better-known arithmetic~$Q$ (e.g., $R$ is not finitely
axiomatizable), it is distinguished by its interpretability properties---see Visser~\cite{vis:r}.

The above-mentioned conditions on theories form an increasing chain. For most of the inclusions in this chain, it is
clear (or at least, reasonably well known) that the inclusions are strict: in particular, there are theories
interpreting~$Q$ that are not sequential (in fact, $Q$ itself is such a theory~\cite{vis:pair}), $R$ does not
interpret~$Q$ (as $Q$ is a finitely axiomatized theory with no finite model, whereas $R$ is locally finitely
satisfiable), and there are essentially undecidable recursively axiomatized theories that do not represent prf.
However, one of these inclusions is not as easy to resolve, leading to the question that motivated this paper:
\begin{Que}\label{que:motiv}
If a theory represents all partial recursive functions, does it interpret Robinson's theory~$R$?
\end{Que}
This may look plausible at first sight: $R$ is a very weak theory that only fixes the values of
elementary arithmetic operations on standard natural numbers, and requires virtually nothing else from the rest of the
model. Now, the definition of representability of prf does provide for natural number constants and definable
functions on them that behave like elementary arithmetic operations as these operations are prf, so everything seems
to be in order.

Despite this, the answer turns out to be negative. The devil is in the ``virtually nothing else'': $R$ does, after
all, involve universally quantified conditions that may look innocuous (in our favourite formulation of~$R$, these
universal quantifiers are bounded by a constant, hence ostensibly ``finite''), but actually turn out to be crucially
important. Using Visser's~\cite{vis:r} characterization, $R$ interprets nontrivial universal theories such as the
theory of infinite discrete linear order. In contrast, prf can be represented in a theory axiomatized purely by
quantifier-free sentences, with no universal quantifiers lurking behind.

We are going to prove that consistent theories with quantifier-free---or even existential---axioms cannot interpret
infinite linear orders and a couple of similar universal theories, and a fortiori, cannot interpret~$R$. This is not
easy to work out directly: the weakness of existential theories---which should intuitively be the reason for
nonexistence of such interpretations---backfires in that we have absolutely no control over the complexity of formulas
that make up potential interpretations, and over the sets they define in models.

Our strategy to solve this problem is to consistently extend the interpreting theory to a theory with quantifier
elimination, using the fact that the empty theory in an arbitrary language~$L$ has a model completion (which we denote
$\ec_L$, being the theory of existentially closed $L$-structures). This fact is well known for relational languages, in
which case $\ec_L$ is the theory of the ``random $L$-structure''. However, we need it for languages with function
symbols, which are mostly neglected in common literature, though the existence of $\ec_L$ was proved in full
generality already by Winkler~\cite{winkler}.

It follows that if a theory is interpretable in a consistent quantifier-free or existential theory, it is weakly
interpretable in $\ec_L$ for some~$L$, and the interpretation can be taken quantifier-free. In order to see that this
is heading in the right direction, we establish a converse result: if an $\exists\forall$ theory is weakly
interpretable in $\ec_L$, it is interpretable in a quantifier-free theory.

We proceed to prove that $\ec_L$ does not, actually, weakly interpret various theories of interest. At this point, we
are heading further and further into model theory, having left the realm of arithmetical theories. It turns out that
our non-interpretability results can be naturally expressed in the language of \emph{classification theory}. Arising
through the work of Shelah~\cite{she:class}, classification theory studies the landscape of ``dividing lines'' between
tame and wild theories, and their structural consequences. Many dividing lines have the following form: a theory is
wild if it has a model that contains a certain complex combinatorial arrangement. Usually, conditions of this form can
be reformulated as (weak) interpretability of a specific $\exists\forall$~theory. For a concrete example, a theory~$T$
has the \emph{strict order property} ($\prpty{SOP}$) if there exists a model $M\model T$, a formula $\fii(\vc x,\vc y)$, and
tuples $\vc a_n\in M$ for $n\in\N$ such that $\fii$ defines in~$M$ a strict partial order, and $\fii(\vc a_n,\vc a_m)$
whenever $n<m$. Otherwise, $T$ is said to have the \emph{no-strict-order property} ($\prpty{NSOP}$).

We observe that theories that can represent recursive functions, as well as consistent extensions of $\ec_L$ for sufficiently rich
languages~$L$, are moderately wild in that they always have the tree property $\prpty{TP_2}$.
However, we will prove
that $\ec_L$ (for arbitrary~$L$) has certain tameness properties: specifically, it has the no-strong-order property
$\prpty{NSOP_3}$ (which implies $\prpty{NSOP}$), and it has elimination of infinity. Using a characterization of
$\prpty{NSOP_1}$ theories by Chernikov and Ramsey~\cite{cher-ram}, we show that it even has the $\prpty{NSOP_1}$
property. On the other hand, theories interpreting~$R$ are firmly on the wild side of all generally considered dividing
lines.

For completeness, the paper also includes discussion of basic model-theoretic properties of~$\ec_L$ in the appendix.

\section{Preliminaries}\label{sec:preliminaries}

Let us first agree on a few bits of general notation. We will use $\N$ and~$\omega$ more or less interchangeably to
denote the set of nonnegative integers; $\N$ may also denote the standard model of arithmetic
$\p{\N,0,\sucf,{+},{\cdot},{<}}$. We denote sequences by angle brackets, and consider them indexed starting from~$0$;
tuples of finite-but-unspecified length will be denoted by placing a bar over a variable name, so that $\vc x$ may
stand for the $n$-tuple $\p{x_0,\dots,x_{n-1}}$.

We will write $F\colon X\pto Y$ to denote that $F$ is a partial function from $X$ to~$Y$. (We use this notation in the
context of partial recursive functions, so virtually always we will have $X=\N^k$, $Y=\N$.)

The notation $t\eqs s$ means that $t$ and~$s$ are syntactically identical terms; we may also apply it to formulas and
other syntactic objects.

\subsection{Theories and interpretations}\label{sec:theor-interpr}

In this paper, a \emph{language} consists of an arbitrary number of relation and function symbols of arbitrary finite
arity (including~$0$: nullary functions are constants, nullary relations are propositional variables; relation and
function symbols of arity $\ge1$ are called \emph{proper}). A \emph{theory}
is a deductively closed set of sentences in a particular language. A theory in language~$L$ is also called an
$L$-theory. We will often consider theories specified by a set
of axioms, in which case the theory is taken to be their deductive closure; we will frequently omit outer universal
quantifiers from axioms. We will generally employ a form of first-order logic that allows \emph{empty models}.

Many considerations in this paper revolve around the notion of \emph{interpretation} of one theory in another, so we
need to be somewhat specific about its meaning. However, since a precise technical definition of interpretations would
get quite lengthy, we advise the reader to consult e.g.\ Visser~\cite[\S2]{vis:r} for the details if necessary; we will
only indicate the main distinctive features below.

Let $T$ be a theory in a language~$L_T$, and $S$ a theory in a language~$L_S$. In its most simple form, a
\emph{translation} $I$ of language $L_T$ into language~$L_S$ is specified by:
\begin{itemize}
\item An $L_S$-formula $\delta_I(x)$ denoting the domain of~$I$.
\item For each relation symbol $R$ of~$L_T$, as well as the equality relation~$=$, an $L_S$-formula $R_I$ of the same
arity.
\item For each function symbol $F$ of~$L_T$ of arity~$k$, an $L_S$-formula $F_I$ of arity~$k+1$.
\end{itemize}
If $\fii$ is an $L_T$-formula, its $I$-translation $\fii^I$ is an $L_S$-formula constructed as follows: we rewrite the
formula in an equivalent way so that function symbols only occur in atomic subformulas of the form
$F(\vc x)=y$, where $x_i,y$ are variables; then we replace each such atomic formula with
$F_I(\vc x,y)$, we replace each atomic formula of the form $R(\vc x)$ with $R_I(\vc x)$,
and we restrict all quantifiers and free variables to objects satisfying~$\delta_I$. We take care to rename bound
variables to avoid variable capture during the process.

A translation $I$ of $L_T$ into~$L_S$ is an \emph{interpretation} of $T$ in~$S$ if $S$ proves:
\begin{itemize}
\item For each function symbol $F$ of~$L_T$, the formula expressing that $F_I$ is total on~$\delta_I$: 
\begin{equation}\label{eq:6}
\forall x_0,\dots,x_{k-1}\,\bigl(\delta_I(x_0)\land\dots\land\delta_I(x_{k-1})\to\exists y\,(\delta_I(y)\land
F_I(\vc x,y))\bigr).
\end{equation}
\item The $I$-translations of all axioms of $T$, and axioms of equality.
\end{itemize}
It follows that $S$ proves the $I$-translations of all sentences provable in~$T$.

The simplified picture of translations and interpretations above actually describes only \emph{one-dimensional},
\emph{parameter-free}, and \emph{one-piece} translations. In the full generality, we allow the following:
\begin{itemize}
\item Translations may be \emph{multi-dimensional}. That is, we use $n$-tuples of $L_S$-objects to represent $L_T$-objects
(where $n$ is a fixed natural number, called the dimension of the translation):
thus, $\delta_I$ has $n$~free variables, $R_I$ has $kn$~free variables for a $k$-ary relation~$R\in L_T$, and
similarly for functions; and when constructing $\fii^I$, each quantifier is replaced with a block of $n$
quantifiers.
\item Translations may use \emph{parameters}. This means that the formulas $\delta_I$, $R_I$, and $F_I$ may include
parameter variables $\vc w$ that are assumed distinct from any proper variables used in the target formulas,
and the specification of $I$ includes an $L_S$-formula $\pi_I(\vc w)$ that describes which parameters are
admissible. Parameters carry through the translation unchanged, so they appear as free variables in~$\fii^I$. The definition
of interpretation is modified so that $S$ proves $\forall\vc w\,(\pi_I(\vc w)\to\fii^I(\vc w))$ for each axiom
$\fii$, and likewise for~\eqref{eq:6}.
\item Translations may be \emph{piece-wise}: the interpreted domain of $L_T$-objects may be stitched together from finitely many
pieces (possibly of different dimensions, and possibly overlapping). Each piece has its own $\delta_I$ formula,
there is a separate $R_I$ formula for each choice of a sequence of pieces for the arguments of~$R$, etc.
\end{itemize}
A translation~$I$ is called \emph{unrelativized} if, on each piece, $\delta_I(\vc x)$ is a tautologically true formula,
and it has \emph{absolute equality} if, on each piece, $\vc x=_I\vc y$ is the formula $\ET_ix_i=y_i$.

Under suitable conditions, we do not need the full generality of interpretations:
\begin{itemize}
\item Assume that $S$ proves the existence of at least two distinct objects. Then whenever $T$ has an interpretation
in~$S$, it also has a one-piece interpretation. (The new interpretation may have larger dimension, but needs no extra
parameters.) This can be achieved by using the pattern of equalities on an extra tuple of variables to distinguish
pieces. For this reason, we will mostly think of interpretations as one-piece, to avoid unnecessary technical baggage.
\item If $T$ has a definable object, then an interpretation of $T$ in~$S$ may be converted to an unrelativized
interpretation by ``equating'' tuples outside the original domain with the definable object. If we do not mind using
extra parameters, the same can be achieved even if $T$ just proves the existence of at least one object. This
construction may not be always desirable, hence relativized interpretations will remain the norm for us.
\item A theory $S$ has (non-functional) \emph{pairing} if there is a formula $\pi(x,y,z)$ such that $S$ proves
\begin{gather*}
\pi(x,y,z)\land\pi(x',y',z)\to x=x'\land y=y'\\
\forall x\,\forall y\,\exists z\,\pi(x,y,z).
\end{gather*}
If $T$ has an interpretation in a theory $S$ with pairing, it also has a one-dimensional interpretation, as we can use
single elements to code tuples.
\end{itemize}

If $I_1$ is a translation of language $L_1$ into~$L_0$, and $I_2$ a translation of language $L_2$ into~$L_1$,
the \emph{composition} $I_1\circ I_0$ is a translation of $L_2$ into~$L_0$, and it is defined in an expected way. Note
that if $I_1$ is an interpretation of a theory $T_1$ in~$T_0$, and $I_2$ an interpretation of $T_2$ in~$T_1$, then
$I_1\circ I_2$ is an interpretation of $T_2$ in~$T_0$.

Let $T$ and~$S$ be theories. Some variants on the notion of interpretation of $T$ in~$S$ are:
\begin{itemize}
\item A \emph{weak interpretation} of $T$ in~$S$ is an interpretation of $T$ in a consistent extension of~$S$ (in the
same language as~$S$), or equivalently, in a completion of $S$.
\item A \emph{cointerpretation} of $T$ in~$S$ is a translation $I$ of language~$L_S$ into~$L_T$ (sic!) such that
$T\vdash\fii^I$ implies $S\vdash\fii$ for every $L_S$-sentence~$\fii$.
\item A \emph{faithful interpretation} of $T$ in~$S$ is an interpretation of $T$ in~$S$ that is at the same time a
cointerpretation of $S$ in~$T$.
\end{itemize}
A theory $T$ is \emph{interpretable} (\emph{weakly interpretable}, \emph{cointerpretable}) in a theory~$S$ if there
exists an interpretation (weak interpretation, cointerpretation, resp.) of $T$ in~$S$.

If $T$ and~$S$ are complete theories, a translation $I$ of $L_T$ in $L_S$ is an interpretation of $T$ in~$S$ iff it is
a weak interpretation iff (assuming $I$ is parameter-free) it is a cointerpretation of $S$ in~$T$.
\begin{Lem}\label{lem:coi-wi}
If $I$ is a weak interpretation of $U$ in~$T$, and $J$ a cointerpretation of $S$ in~$T$, then $J\circ I$ is a weak
interpretation of $U$ in~$S$.
\noproof\end{Lem}

An interpretation $I$ of $T$ in~$S$, as defined, is a syntactic transformation of formulas provable in~$T$ into
formulas provable in~$S$. However, it can be also viewed semantically: it provides a uniform way of building
``internally definable'' models of~$T$ out of models of~$S$.

Assume first $I$ is a parameter-free one-piece interpretation with absolute equality, and let $M\model S$. We construct
a model $M^I\model T$ as follows: if $I$ is $n$-dimensional, the domain of $M^I$ is $\delta_I(M)=\{\vc a\in
M^n:M\model\delta_I(\vc a)\}$; a $k$-ary relation symbol $R\in L_T$ is realized in~$M^I$ by $\{\p{\vc a_0,\dots,\vc
a_{k-1}}\in\delta_I(M)^k:M\model R_I(\vc a_0,\dots,\vc a_{k-1})\}$, and similarly, a $k$-ary function symbol $F\in
L_T$ is realized by the function whose graph is the subset of $(\delta_I(M))^{k+1}$ defined in~$M$ by the formula~$F_I$.

Next, if $I$ does not have absolute equality, we build the structure as before, and let $M^I$ be its quotient by the
binary relation defined on it by the formula $=_I$; this relation is in fact a congruence, as $S$ proves the
translations of equality axioms.

If $I$ is a piece-wise interpretation, we construct the domain of $M^I$ as the disjoint union of the finitely many
pieces, each defined as above; we define relations and functions in the appropriate way.

Finally, if $I$ is an interpretation with parameters, we will not obtain a single model~$M^I$, but one model for each
choice of parameters: that is, if $\vc a$ is a tuple such that $M\model\pi_I(\vc a)$, then $M^{I,\vc a}$ is a model
of~$T$ built from the expanded structure $\p{M,\vc a}$ by the procedure above.

\subsection{Representation of recursive functions}\label{sec:repr-recurs-funct}

The notion of representable\footnote{In the terminology of~\cite{tmr}, \emph{definable}. We reserve the latter word for
something else, in accordance with current standard usage.} predicates and
functions in first-order theories was introduced in~\cite{tmr}. We summarize it below, with a few inessential
modifications. (Warning: we are going to relax the definition a bit later in this section.)
\begin{Def}\label{def:rep-rec}
Let $T$ be a theory in a language~$L$, and $\sigma=\{\num n:n\in\N\}$ a fixed \emph{sequence of numerals}: i.e., a
sequence of closed terms~$\num n$ such that
\[T\vdash\num n\ne\num m\]
for $n,m\in\N$, $n\ne m$.

A recursive predicate (rp) $P\sset\N^k$ is \emph{represented} in~$T$ w.r.t.~$\sigma$ by a formula $\fii(x_0,\dots,x_{k-1})$ if
\begin{align*}
\p{n_0,\dots,n_{k-1}}\in P&\implies T\vdash\phantom\neg\fii(\num{n_0},\dots,\num{n_{k-1}}),\\
\p{n_0,\dots,n_{k-1}}\notin P&\implies T\vdash\neg\fii(\num{n_0},\dots,\num{n_{k-1}})
\end{align*}
for all $n_0,\dots,n_{k-1}\in\N$.

A partial recursive function (prf) $F\colon\N^k\pto\N$ is \emph{represented} w.r.t.~$\sigma$ by a formula
$\fii(\vc x,y)$ if 
\[T\vdash\fii(\num{n_0},\dots,\num{n_{k-1}},y)\eq y=\num m\]
whenever $n_0,\dots,n_{k-1},m\in\N$ are such that $F(\vc n)=m$.

A set $\mathcal R$ of prf and rp is representable in~$T$ if there exists a
sequence of numerals~$\sigma$ such that each member of~$\mathcal R$ is representable in~$T$ w.r.t.~$\sigma$.
\end{Def}

In fact, \cite{tmr} only considers representation of total recursive functions (trf), but it can be obviously generalized to partial
functions in the indicated fashion. Likewise, we can generalize representation of rp to
representation of \emph{disjoint pairs of r.e.\ predicates (dprp)}: such a disjoint pair $\p{P^+,P^-}$, where
$P^+,P^-\sset\N^k$, is represented by a formula~$\fii(\vc x)$ if
\begin{align*}
\p{n_0,\dots,n_{k-1}}\in P^+&\implies T\vdash\phantom\neg\fii(\num{n_0},\dots,\num{n_{k-1}}),\\
\p{n_0,\dots,n_{k-1}}\in P^-&\implies T\vdash\neg\fii(\num{n_0},\dots,\num{n_{k-1}})
\end{align*}
for all $n_0,\dots,n_{k-1}\in\N$. We identify any relation $P\sset\N^k$ with the disjoint pair $\p{P,\N^k\bez P}$.

Notice that a representation of a rp $P\sset\N^k$ is essentially the same as a representation of its
characteristic function $\chi_P\colon\N^k\to\N$; likewise for disjoint pairs (their characteristic functions are
partial). Consequently, representability of all prf in~$T$ implies representability of all trf and representability of all dprp; in
turn, either of the latter two properties implies representability of all rp.

The definition of representation of functions does not demand anything from $\fii(\vc x,y)$ when $\vc x$ is not one
of the tuples $\vc{\num n}$ in the domain of the original function. However, if $\fii(\vc x,y)$ represents a partial
function $F\colon\N^k\pto\N$ in~$T$, we may define
\begin{align*}
\fii'(\vc x,y)&\eqs\forall z\,(\fii(\vc x,z)\eq z=y),\\
\fii''(\vc x,y)&\eqs\fii'(\vc x,y)\lor\bigl(y=\num0\land\neg\exists z\,\fii'(\vc x,z)\bigr).
\end{align*}
Then $\fii'$ and $\fii''$ also represent $F$ in~$T$; moreover, $\fii'$ is $T$-provably a partial function, and $\fii''$
is $T$-provably a total function. Thus, we could have included either condition in the definition with no ill effects.

A desirable condition that we did not include in the definition is that the sequence of numerals~$\sigma$ be
\emph{recursive}:
that is, we can compute the term $\num n$ on input~$n$. For most purposes, this is actually redundant if $T$ can
represent recursive functions with respect to~$\sigma$: using a formula representing the (recursive) successor
function $\sucf(n)=n+1$, we can build a recursive sequence of formulas~$\fii_n(x)$ that define~$\num n$.

Definition~\ref{def:rep-rec} formally makes sense for representation of \emph{arbitrary} predicates or partial functions in~$T$.
However, there is little point in that: if $T$ is recursively axiomatizable, and the given numeral sequence is
recursive (or if we can represent~$\sucf$), then all predicates and total functions represented
in~$T$ are actually recursive, and each partial function represented in~$T$ extends to a partial recursive
function represented in~$T$. (This is not necessarily true for non-recursive numeral sequences, see
Proposition~\ref{prop:nonrecnumseq}.)

The primary reason for discussing representability of recursive functions in~\cite{tmr} is that it implies essential
undecidability. We include the argument below for completeness.
\begin{Prop}\label{prop:rep-ess-undec}
If the set of all unary rp is representable in a theory~$T$ w.r.t.\ a recursive sequence of numerals,
then $T$ is essentially undecidable.
\end{Prop}
\begin{Pf}
Let $S\Sset T$ be decidable. This makes the predicate
\[P(n)\iff\text{$n$ is the G\"odel number of a formula $\fii(x)$ s.t.\ $S\vdash\neg\fii(\num n)$}\]
recursive, hence $P$ is represented in~$S$ by a formula $\fii(x)$. Let $n=\code\fii$ be its G\"odel number. If $\neg
P(n)$, then $S\vdash\neg\fii(\num n)$ by representability, hence $P(n)$ by the definition of~$P$, which is a
contradiction. Thus, $P(n)$. Then $S\vdash\neg\fii(\num n)$ by the definition of~$P$, and $S\vdash\fii(\num n)$ by
representability, hence $S$ is inconsistent.
\end{Pf}

Again, the assumption of recursivity of the numeral sequence in Proposition~\ref{prop:rep-ess-undec} may be replaced
with representability of~$\sucf$. However, it cannot be dropped entirely, as shown in the appendix
(Propositions \ref{prop:nonrecnumseq} and~\ref{prop:nonrec-all}).

Likewise, it is essential in Proposition~\ref{prop:rep-ess-undec} that all unary rp are representable at once:
we show in Proposition~\ref{prop:unif-rec-dec} that any finite (or uniformly recursive) set of rp and trf is representable in a
decidable theory. In contrast, there is one fixed unary dprp (or: prf) whose representability in a theory w.r.t.\ a
recursive numeral sequence implies essential undecidability: in fact, any recursively inseparable pair has this
property.

The reader may have realized that representation of recursive functions and predicates in~$T$ amounts to an
interpretation of a particular theory in~$T$. We now make this connection explicit.
\begin{Def}\label{def:tprf}
Let $\mathcal R$ be a set of prf and dprp. The language $L_{\mathcal R}$ consists of constants $\{\num n:n\in\N\}$,
function symbols $\num F$ of appropriate arity for every prf $F\in\mathcal R$, and likewise relation symbols $\num P$
for every dprp $P\in\mathcal R$. The theory
$\trep{\mathcal R}$ in language $L_{\mathcal R}$ is axiomatized by
\[\num n\ne\num m\]
for $n\ne m\in\N$;
\[\num F(\num{n_0},\dots,\num{n_{k-1}})=\num m\]
for each $k$-ary function $F\in\mathcal R$, and $n_0,\dots,n_{k-1},m\in\N$ such that $F(\vc n)=m$; and for each
$k$-ary disjoint pair $P=\p{P^+,P^-}\in\mathcal R$, the axioms
\[\num P(\num{n_0},\dots,\num{n_{k-1}})\]
for $\p{n_0,\dots,n_{k-1}}\in P^+$, and
\[\neg\num P(\num{n_0},\dots,\num{n_{k-1}})\]
for $\p{n_0,\dots,n_{k-1}}\in P^-$.
This definition also applies to rp $P\sset\N^k$ using their identification with dprp $\p{P,\N^k\bez P}$.

Note that that the theory $\trep{\mathcal R}$ is axiomatized by open (= quantifier-free) sentences. 

Let $\sprf$, $\strf$, $\sdprp$, and $\srp$ denote the sets of all prf, trf, dprp, and rp, respectively (where we
consider $\strf\sset\sprf$ and $\srp\sset\sdprp$).
Since $\trep\sdprp$ is
included in an extension of $\tprf$ by quantifier-free definitions, we will use $\tprf$ as a proxy for
$\trep{\sprf\cup\sdprp}$. 

For convenience, we also consider a finite-language formulation of~$\tprf$. Let $U(x,y)$ be the prf defined by
\begin{align*}
U(0,m)&=m+1,\\
U(n+1,m+1)&=[n,m],\\
U([n,m]+1,0)&\simeq\varphi_n(m),
\end{align*}
where $[n,m]$ denotes a recursive bijective pairing function $\N^2\to\N$ (e.g., the Cantor pairing function $(n+m)(n+m+1)/2+n$),
and $\varphi_n(m)$ a partial recursive numbering of unary prf. Let $\tprfu$ be the fragment of~$\tprf$ in the language
$\p{\num0,\num U}$; it can be axiomatized by
\begin{gather*}
\num0\ne S(\num0),\\
\num U(S^n(\num0),S^m(\num0))=S^k(\num0)
\end{gather*}
for all $n,m,k\in\N$ such that $U(n,m)=k$, where $S(x)$ denotes $\num U(\num0,x)$.
\end{Def}
\begin{Lem}\label{lem:tprfu}
$\tprf$ is included in an extension of~$\tprfu$ by definitions of function symbols by terms, thus a theory interprets
$\tprf$ iff it interprets~$\tprfu$.
\end{Lem}
\begin{Pf}
We can read the definition of $U$ backwards to obtain definitions of $S$, $[x,y]$, and $\varphi_n(x)$ in terms of $U$
and~$0$: $S(x)=U(0,x)$, $[x,y]=U(S(x),S(y))$, $\varphi_n(x)=U(S([S^n(0),x]),0)$. Then any prf $F\colon\N^k\pto\N$ can be
written in the form $F(x_0,\dots,x_{k-1})=\varphi_n([x_0,[x_1,\cdots[x_{k-2},x_{k-1}]\cdots]])$ for a suitable~$n$.
\end{Pf}

Using the above-mentioned fact that representations of (partial) functions may be assumed to be actual definable
functions, we see:
\begin{Obs}\label{obs:rep-int}
A set $\mathcal R$ of prf and dprp is representable in a theory~$T$ according to
Definition~\ref{def:rep-rec} iff $\trep{\mathcal R}$ is interpretable in~$T$ by a one-piece one-dimensional parameter-free
interpretation~$I$ with absolute equality such that each $\num n^I$ is definable in~$T$ by a closed term.
\noproof\end{Obs}

Now, the restrictions on the interpretation in Observation~\ref{obs:rep-int} are mostly irrelevant and arbitrary; as we are
looking at the concept of representations from the viewpoint of interpretability, it seems we obtain a cleaner concept
if we just drop them:
\begin{Def}\label{def:rep-int}
A \emph{loose representation} of $\mathcal R\sset\sprf\cup\sdprp$ in a theory~$T$ is an interpretation of
$\trep{\mathcal R}$ in~$T$.
\end{Def}

In particular, a theory $T$ loosely represents all prf iff it interprets the theory~$\tprfu$.

\subsection{The theory $R$}\label{sec:theory-r}

Robinson's theory $R$ was originally defined in~\cite{tmr}. Some inessential variants (mutually interpretable) of the
theory appear in the literature; we prefer the following form in this paper.
\begin{Def}\label{def:r}
Let $R$ denote the theory in the language $L_R=\p{0,\sucf,+,\cdot,<}$ axiomatized by
\begin{gather}
\label{eq:7}  \num n+\num m=\num{n+m},\\
\label{eq:11} \num n\cdot\num m=\num{nm},\\
\label{eq:12} x<\num n\eq x=\num0\lor\dots\lor x=\num{n-1}
\end{gather}
for all~$n,m\in\N$, where $\num n\eqs\sucf^n(0)$.
\end{Def}
(In particular, note that axiom~\eqref{eq:12} for $n=0$ states $\neg(x<0)$.) It is easy to show that $R$ implies $\num
n\ne\num m$ for distinct $n,m\in\N$.

Observe that an $L_R$-structure is a model of~$R$ iff it contains the standard model~$\N$ as an initial (i.e., closed
downward under~$<$) substructure.

As usual, \emph{bounded quantifiers} are introduced in~$L_R$ as the short-hands
\begin{align*}
\exists y<t(\vc x)\,\fii(\vc x,y)&\eqs\exists y\,(y<t(\vc x)\land\fii(\vc x,y)),\\
\forall y<t(\vc x)\,\fii(\vc x,y)&\eqs\forall y\,(y<t(\vc x)\to\fii(\vc x,y)),
\end{align*}
where $t$ is a term not containing the variable~$y$. An $L_R$-formula $\fii(\vc x)$ is $\Delta_0$ (or \emph{bounded})
if all quantifiers in~$\fii$ are bounded. A formula is $\Sigma_1$ if it consists of a block of existential quantifiers
followed by a $\Delta_0$~formula.
\begin{Prop}\label{prop:r-sig1}
$R$ proves all $\Sigma_1$~sentences true in the standard model~$\N$. Conversely, it can be axiomatized by a set
of true (universal) $\Delta_0$~sentences.
\noproof\end{Prop}

As already proved in~\cite{tmr} (for the original, slightly stronger definition of the theory), $R$ can represent
recursive functions. We briefly sketch the argument below for completeness.
\begin{Prop}\label{prop:r-prf}
Every prf $F\colon\N^k\pto\N$ is representable in~$R$ by a $\Sigma_1$~formula w.r.t.\
the usual sequence of numerals as in Definition~\ref{def:r}.
\end{Prop}
\begin{Pf}
The graph $\{\p{\vc x,y}:F(\vc x)=y\}$ is definable in~$\N$ by a $\Sigma_1$~formula of the form $\exists z\,\tet(\vc
x,y,z)$, where $\tet\in\Delta_0$. Put
\begin{align*}
\alpha(w)&\eqs0<w\land\forall z<w\,(\sucf(z)=w\lor\sucf(z)<w),\\
\fii(\vc x,y)&\eqs\exists w,z\,\bigl(\alpha(w)\land y<w\land z<w\land\tet(\vc x,y,z)
\land\forall y',z'<w\,(\tet(\vc x,y',z')\to y'=y)\bigr).
\end{align*}
One can check
\begin{equation}\label{eq:13}
R\vdash\alpha(w)\to w=\num1\lor\dots\lor w=\num r\lor\num r<w
\end{equation}
for any~$r\in\N$.

We claim that $\fii$ represents $F$ in~$R$. Assume $F(\vc{\num n})=\num m$. On the one hand, $\fii(\vc{\num n},\num m)$
is a true $\Sigma_1$~sentence, and as such it is provable in~$R$. On the other hand, fix $r\in\N$ that witnesses
the $z$ quantifier in $\N\model\fii(\vc{\num n},\num m)$. Working in~$R$, assume $\fii(\vc{\num n},y')$, we need
to show $y'=\num m$. Let $w',z'$ witness the existential quantifiers in $\fii(\vc{\num n},y')$. Using~\eqref{eq:13},
either $w'$ equals a standard numeral, or $w'>\num m,\num r$. In the latter case, $\theta(\vc{\num n},\num m,\num r)$
implies $y'=\num m$ as needed. In the former case, $y',z'<w'$ are also standard. It again follows that $y'=\num m$, as
otherwise $\neg\theta(\vc{\num n},y',z')$ would be a true $\Delta_0$ sentence, thus provable in~$R$.
\end{Pf}
Consequently, $R$ is essentially undecidable.

It is easy to see that $R$ (therefore any theory interpretable in~$R$) is \emph{locally finitely satisfiable}%
\footnote{This terminology from \cite{vis:r} is unrelated to the notion of a type being finitely satisfiable.}%
, i.e., every finite subset has a finite model:
indeed, if we identify all elements of~$\N$ above $b+1$, we obtain a model satisfying \eqref{eq:7}, \eqref{eq:11}, and
\eqref{eq:12} for $n\le b$. Visser~\cite{vis:r} proved a striking converse to this observation:
\begin{Thm}\label{thm:r-locfinsat}
Every locally finitely satisfiable, recursively axiomatizable theory in a finite language is interpretable in~$R$,
using a one-piece one-dimensional parameter-free interpretation.
\noproof\end{Thm}

Since relational $\exists\forall$~sentences have the finite model property, this in particular implies
that $R$ interprets any consistent theory axiomatized by a recursive set of $\exists\forall$~sentences in a finite
relational language.

\subsection{Model theory}\label{sec:model-theory}

Since this paper is intended to be accessible to a non-model-theoretic audience (and the author is not a model theorist
either), it will only assume modest prerequisites in model theory---mostly common knowledge among logicians. We will
review a few selected topics in more detail below; the material needed should be covered by a textbook such
as~\cite{hodges:sh}, except that we will also need a few concepts from classification theory, which will be explained
in the next section.

First, let us start with a few basic conventions. Recall that we allow models to be empty, and that we denote finite
tuples as~$\vc x$. For any structure~$M$, we denote by $\diag(M)$ its diagram: the set of quantifier-free sentences true
in~$M$ in the language of~$M$ augmented with constants for each element of~$M$. By a slight abuse of language, we will
also use this notation to denote the set of quantifier-free sentences true in~$M$ in its \emph{original} language, if
every element of~$M$ is the value of a closed term (i.e., if $M$ is $0$-generated).

Even though we normally work with one-sorted logic, the following construction is best thought of as yielding a
multi-sorted structure. For any structure~$M$, let $M^\teq$ be the structure that has $M$ itself as one of its sorts,
and for each equivalence relation $E(\vc x,\vc y)$ on~$M^n$ definable without parameters in~$M$, it has a sort whose elements are
the equivalence classes of~$E$; the structure includes the projection function to this sort from~$M^n$. It is easy to
see that each such equivalence relation is definable in~$M$ by a formula that \emph{provably} defines an equivalence
relation in predicate logic; thus, the following makes sense: for any theory~$T$, let $T^\teq$ be the multi-sorted
theory whose models are exactly the structures $M^\teq$ for $M\model T$. (Officially, $M^\teq$ and $T^\teq$ can be
coded in a suitable one-sorted language.) Note that $T^\teq$ is interpretable in~$T$, and any interpretation of another
theory $S$ in~$T$ can be made into an interpretation with absolute equality of $S$ in~$T^\teq$.

Since we will work a lot with model completions, let us recall the related background. Let $\mathcal K$ be a class of
structures in the same language. A model $M\in\mathcal K$ is
\emph{existentially closed (e.c.) in~$\mathcal K$} if for every model $N\Sset M$ such that $N\in\mathcal K$, we
have $M\preceq_1N$: i.e., every existential formula with parameters from~$M$ which is satisfied in~$N$ is already
satisfied in~$M$. We will often speak of (absolutely) e.c.\ models without reference to~$\mathcal K$, in which case it
is understood that $\mathcal K$ is the class of all models in the given language. An e.c.\ model of a theory~$T$ is an
e.c.\ structure in the class of models of~$T$. If $T$ is a $\forall\exists$-axiomatized theory, then every model
$M\model T$ embeds in an e.c.\ model of~$T$. (More generally, this holds for any class $\mathcal K$ closed under
limits of chains.)

A theory~$T$ is \emph{model-complete} if all models $M\model T$ are e.c.\ models of~$T$; this implies the stronger
condition that for all $M,N\model T$, $M\sset N$ implies $M\preceq N$. Equivalently, $T$ is model-complete iff every
formula~$\fii$ is in~$T$ equivalent to an existential formula; it is enough to test this for universal formulas~$\fii$.
A stronger condition is that $T$ has \emph{quantifier elimination}, meaning that every formula~$\fii$ is in~$T$
equivalent to a quantifier-free formula; it is enough to test this for existential formulas~$\fii$ with only one
quantifier. Any model-complete theory~$T$ is axiomatizable by $\forall\exists$~sentences.

Theories $T$ and~$S$ in the same language are \emph{companions} if every model of $T$ embeds in a model of~$S$, and
vice versa; equivalently, $T_\forall=S_\forall$, where $T_\forall$ denotes the universal fragment of~$T$. A \emph{model
companion} of a theory~$T$ is a model-complete theory~$T^*$ that is a companion of~$T$. There are theories with no
model companion (e.g., the theory of groups), but if a theory~$T$ has a model companion~$T^*$, it is unique: the models
of~$T^*$ are exactly the e.c.\ models of~$T_\forall$. A theory has a model companion iff the class of e.c.\ models
of~$T_\forall$ is elementary. Notice that a model companion of~$T$ is the same thing as a model companion of~$T_\forall$,
hence we can as well restrict attention to universal theories~$T$.

A \emph{model completion} of a theory~$T$ is a model companion $T^*$ of~$T$ such that for every $M\model T$, the theory
$T^*+\diag(M)$ is complete. Equivalently, a model companion $T^*$ of~$T$ is a model completion of~$T$ iff $T$ has the
amalgamation property (cf.\ Definition~\ref{def:fra}). If $T$ is a universal theory (which is the case we are primarily
interested in), a companion $T^*$ of~$T$ is a model completion of~$T$ iff $T^*$ has quantifier elimination.

A convenient trick when studying models of a complete theory~$T$ is to use \emph{monster models}. A monster model
of~$T$ is a model $\monster\model T$ sufficiently rich so that all models we need to discuss can be assumed to be
submodels of~$\monster$; in order for this to work, we make $\monster$ highly saturated: to be specific, let us posit
that $\monster$ is $\kappa$-saturated (i.e., every type over $<\kappa$ parameters from~$\monster$ is realized
in~$\monster$) and strongly $\kappa$-homogeneous (i.e., every partial elementary self-map of~$\monster$ of size
$<\kappa$ extends to an automorphism of~$\monster$), where $\kappa$ is a ``large'' cardinal number (in particular,
larger than the size of the language, as well as any models of~$T$ that we are going to encounter during the argument).
This also implies that $\monster$ is $\kappa^+$-universal (every model of~$T$ of size $\le\kappa$ elementarily embeds
in~$\monster$). (If it were not for foundational issues that we prefer not to be dragged into, we could even take
$\monster$ as an ``$\mathrm{Ord}$-saturated'' model: a proper class model of~$T$ saturated w.r.t.\ types over any
\emph{set} of parameters.) Having fixed the monster model~$\monster$, a \emph{small set} is a subset of~$\monster$ of
size $<\kappa$ (likewise for sequences and other similar objects); a \emph{small model} is an elementary submodel
of~$\monster$ of size $<\kappa$.

\subsection{Classification theory}\label{sec:class-theory}

Stability theory and the more general
classification theory was initially developed by Shelah~\cite{she:stab,she:class} (with some notions pioneered by
Morley~\cite{morley}); one of its main themes is identifying useful ``dividing lines'' between tame and wild theories.
The dividing lines we are going to mention here are mostly combinatorial properties based on the appearance of certain
arrangements of points and definable sets in models; for other kind of dividing lines (variants of stability based on
counting of types), see Appendix~\ref{sec:dependence-language}.

While model theorists prefer to work with complete theories, the properties below are all stated in such a way that a
theory~$T$ has a ``tameness'' property~$P$ iff every completion of~$T$ has property~$P$. Also, it will be generally the
case that $T$ has a (tameness) property~$P$ iff every countable-language fragment of~$T$ has property~$P$.

\begin{figure}
\centering
\includegraphics{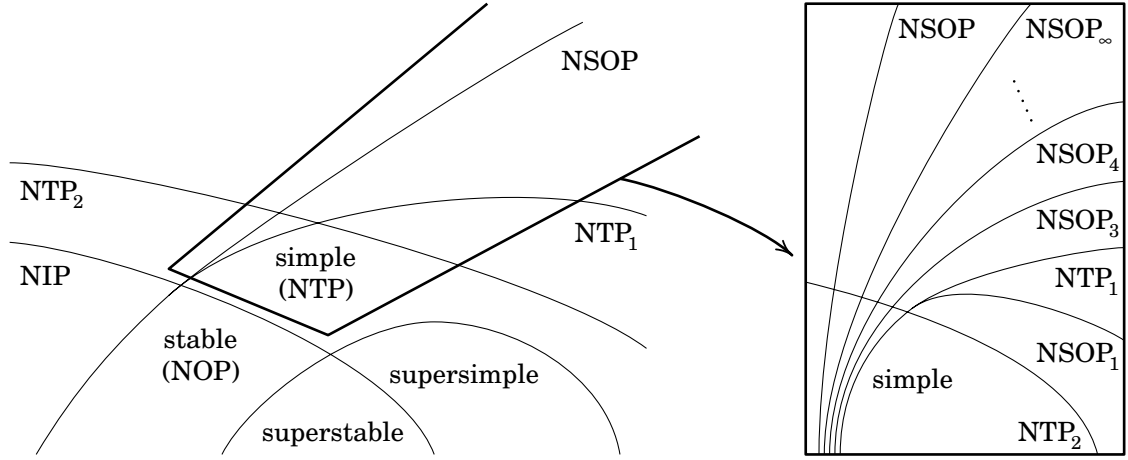}
\caption{Main dividing lines}
\label{fig:main}
\end{figure}
For an overview of inclusions among the properties below, see Figure~\ref{fig:main}.

A theory $T$ has the \emph{order property} ($\prpty{OP}$) if there exists
a formula $\fii(\vc x,\vc y)$ (where $\vc x$ and $\vc y$ are tuples of the same length), a model $M\model T$, and a
sequence of tuples $\{\vc a_i:i\in\omega\}$ in~$M$ such that
\[M\model\fii(\vc a_i,\vc a_j)\iff i<j\]
for all $i,j\in\omega$; otherwise, $T$ has the \emph{no-order property} ($\prpty{NOP}$). It turns out that $T$ has
$\prpty{NOP}$ if and only if it is \emph{stable} (see Appendix~\ref{sec:dependence-language}).

Beware of the terminological peculiarity that the base form of this condition on theories ($\prpty{OP}$) is ``negative''
(wild), whereas the corresponding ``positive'' (tame) condition is denoted as its negation ($\prpty{NOP}$). All
properties below follow the same naming pattern.

A theory~$T$ has the \emph{independence property} ($\prpty{IP}$) if there is a formula $\fii(\vc x,\vc y)$, a
model~$M\model T$, and tuples $\{\vc a_i:i\in\omega\}$ and $\{\vc b_I:I\sset\omega\}$ in~$M$ such that
\[M\model\fii(\vc a_i,\vc b_I)\iff i\in I\]
for all $i\in\omega$ and $I\sset\omega$. Otherwise, $T$ is $\prpty{NIP}$ (also called \emph{dependent}).

A theory~$T$ has the \emph{strict order property} ($\prpty{SOP}$) if there is a formula $\fii(\vc x,\vc y)$, a
model~$M\model T$, and tuples $\{\vc a_i:i\in\omega\}$ in~$M$ such that
\[M\model\exists\vc x\,\bigl(\fii(\vc x,\vc a_i)\land\neg\fii(\vc x,\vc a_j)\bigr)\iff i<j\]
for all $i,j\in\omega$; equivalently, $T$ is $\prpty{SOP}$ iff there is a formula $\fii(\vc x,\vc y)$ that
$T$-provably defines a strict partial order, and there is a model $M\model T$ in which the partial order defined
by~$\fii$ has an infinite chain. Otherwise, $T$ is $\prpty{NSOP}$.

A theory is $\prpty{NOP}$ (stable) if and only if it is both $\prpty{NIP}$ and $\prpty{NSOP}$.

Recall that $X^{<\omega}$ denotes the set of finite sequences with entries from~$X$, ordered by the initial subsequence
relation (which we write as $t\sset s$) to form an
$X$-branching tree; $X^\omega$ is the corresponding set of infinite sequences (which are branches of the tree). If
$s\in X^{<\omega}\cup X^\omega$, and $n<\Lh(s)$, then $s\res n$ is the initial subsequence of $s$ of length~$n$;
if $s\in X^{<\omega}$ and $x\in X$, then $s\cat x$ is $s$ extended with a new entry $x$ at the end. For clarity, we will
write von Neumann numerals as $\jvn n=\{0,\dots,n-1\}$.

A set of formulas is \emph{$k$-inconsistent} if each $k$-element subset is inconsistent.

A theory~$T$ has the \emph{tree property} ($\prpty{TP}$) if there is a formula $\fii(\vc x,\vc y)$, a
model~$M\model T$, tuples $\{\vc a_s:s\in\omega^{<\omega}\}$ in~$M$, and $k\ge2$ such that
\begin{itemize}
\item for each $\sigma\in\omega^\omega$, the type $\{\fii(\vc x,\vc a_{\sigma\res n}):n<\omega\}$ is consistent, and
\item for each $s\in\omega^{<\omega}$, $\{\fii(\vc x,\vc a_{s\cat i}):i<\omega\}$ is $k$-inconsistent.
\end{itemize}
Otherwise, $T$ is called $\prpty{NTP}$ or \emph{simple}. Simplicity can also be equivalently defined in terms of
properties of forking; there is a related stronger condition called \emph{supersimplicity}, see e.g.\
\cite{wagner}. Stable theories are simple, and simple theories are $\prpty{NSOP}$.

The tree property has two important variants. A theory~$T$ has the \emph{tree property $\prpty{TP_1}$} if there is a
formula $\fii(\vc x,\vc y)$, a model~$M\model T$, and tuples $\{\vc a_s:s\in\omega^{<\omega}\}$ in~$M$ such that
\begin{itemize}
\item for each $\sigma\in\omega^\omega$, $\{\fii(\vc x,\vc a_{\sigma\res n}):n<\omega\}$ is consistent, and
\item for each incomparable $s,t\in\omega^{<\omega}$, $\{\fii(\vc x,\vc a_s),\fii(\vc x,\vc a_t)\}$ is inconsistent.
\end{itemize}
$T$ has the \emph{tree property $\prpty{TP_2}$} if there is a
formula $\fii(\vc x,\vc y)$, a model~$M\model T$, and tuples $\{\vc a_{n,i}:n,i\in\omega\}$ in~$M$ such that
\begin{itemize}
\item for each $\sigma\in\omega^\omega$, $\{\fii(\vc x,\vc a_{n,\sigma(n)}):n<\omega\}$ is consistent, and
\item for each $n,i,j\in\omega$ such that $i<j$, $\{\fii(\vc x,\vc a_{n,i}),\fii(\vc x,\vc a_{n,j})\}$ is
inconsistent.
\end{itemize}
As usual, if $T$ is not $\prpty{TP}_i$, it is $\prpty{NTP}_i$. A theory is $\prpty{NTP}$ if and only if it is both
$\prpty{NTP_1}$ and $\prpty{NTP_2}$. All $\prpty{NIP}$ theories are $\prpty{NTP_2}$, and all $\prpty{NTP_1}$ theories
are $\prpty{NSOP}$.

The region between simple and $\prpty{NSOP}$ theories is further stratified by levels of the \emph{strong order
property}. For $k\ge3$, a theory $T$ has the \emph{strong order property $\prpty{SOP}_k$} if there is a formula
$\fii(\vc x,\vc y)$, a model $M\model T$, and tuples $\{\vc a_i:i<\omega\}$ in~$M$ such that $M\model\fii(\vc
a_i,\vc a_j)$ for all $i<j<\omega$, but
\begin{equation}\label{eq:14}
\{\fii(\vc x_0,\vc x_1),\fii(\vc x_1,\vc x_2),\dots,\fii(\vc x_{k-2},\vc x_{k-1}),\fii(\vc x_{k-1},\vc x_0)\}
\end{equation}
is inconsistent; otherwise, $T$ has $\prpty{NSOP}_k$. A theory $T$ has the \emph{strong order property
$\prpty{SOP_\infty}$} if there are data as above such that \eqref{eq:14}
is inconsistent for all $k\ge3$; otherwise, $T$ has $\prpty{NSOP_\infty}$. For any theory $T$, we have
\[\prpty{NTP_1}\implies\prpty{NSOP_3}\implies\prpty{NSOP_4}\implies\cdots\implies\prpty{NSOP_\infty}\implies\prpty{NSOP}.\]
We warn the reader that usage of the abbreviations $\prpty{(N)SOP_\infty}$ and $\prpty{(N)SOP}$ varies in the literature.

Notice that the definition of $\prpty{NSOP}_k$ above is only interesting for $k\ge3$, as stated: taking it blindly
for $k=2$ would give a condition equivalent to $\prpty{NOP}$, and for $k=1$ a condition false for any theory with
infinite models. Instead, the names $\prpty{NSOP_1}$ and $\prpty{NSOP_2}$ were given ad~hoc to variants of the tree
property that fit nicely in the picture. Since $\prpty{NSOP_2}$ is equivalent to $\prpty{NTP_1}$, we will not
bother to define it separately. A theory~$T$ has property $\prpty{SOP_1}$ if there is a formula $\fii(\vc x,\vc y)$, a model~$M\model T$, and tuples $\{\vc a_s:s\in\jvn2^{<\omega}\}$ in~$M$ such that
\begin{itemize}
\item for each $\sigma\in\jvn2^\omega$, $\{\fii(\vc x,\vc a_{\sigma\res n}):n<\omega\}$ is consistent, and
\item for each $s,t\in\jvn2^{<\omega}$, if $s\cat0\sset t$, then $\{\fii(\vc x,\vc a_{s\cat1}),\fii(\vc x,\vc a_t)\}$ is inconsistent;
\end{itemize}
otherwise, $T$ is $\prpty{NSOP_1}$. We have
\[\prpty{NTP}\implies\prpty{NSOP_1}\implies\prpty{NSOP_2}\iff\prpty{NTP_1}\implies\prpty{NSOP_3}\implies\cdots\]
for any theory~$T$.

We observe that each of the combinatorial properties above ($\prpty{NOP}$, $\prpty{NSOP}$, $\prpty{NIP}$, $\prpty{NTP}_{1,2}$, $\prpty{NSOP}_k$, $\prpty{NSOP}_\infty$) can be expressed as weak non-interpretability of a
particular recursively $\exists\forall$-axiomatized theory in a finite relational language. (Incidentally, notice that
any such theory is interpretable in~$R$ by Theorem~\ref{thm:r-locfinsat}, hence $R$ is ``wild'' according to all of these
dividing lines.)

For example, by compactness, a
theory~$T$ has $\prpty{SOP}$ iff it has a model with a definable strict order (on $k$-tuples, for some~$k$) with
arbitrarily long finite chains. It makes no difference if the defining formula is allowed extra parameters, or if we
allow to relativize the domain of the order. Thus, $T$ has $\prpty{SOP}$ iff it has a completion that interprets the
theory $T_\prpty{SOP}$ of strict orders with arbitrarily long chains, axiomatized by
\begin{gather*}
\forall x,y,z\,(x<y\land y<z\to x<z),\\
\forall x\,\neg(x<x),\\
\exists x_0,\dots,x_n\,\ET_{i<n}x_i<x_{i+1}
\end{gather*}
for $n\in\omega$.

For a more complicated example, $T$ has $\prpty{TP_1}$ iff it weakly interprets the theory in a language with
a single binary relation $R(x,y)$, and axioms
\[
\exists\{y_s:s\in\jvn n^{\le n}\}\,\Bigl(\ET_{s\in\jvn n^n}\exists x\,\ET_{i\le n}R(x,y_{s\res i})
  \land\ET_{\substack{s,t\in\jvn n^{\le n}\\s\nsset t\nsset s}}\forall x\,\neg\bigl(R(x,y_s)\land R(x,y_t)\bigr)\Bigr)
\]
for $n\in\omega$.

\section{Model completion of the empty theory}\label{sec:model-compl-empty}

Recall that our original motivation was to find a theory~$T$ that represents prf, but does not interpret~$R$. Now, the
weaker $T$ is, the lower its chances of interpreting~$R$, so the obvious choice is to take $T=\tprf$. This theory,
axiomatized by quantifier-free sentences, essentially just states that the universe includes a copy of a certain model
based on the integers; it does not take a big leap of faith to surmise it is too weak to interpret much of anything. It
is, however, another matter to actually prove this. A possible strategy is to consider an arbitrary translation~$I$ of
the language of~$R$ into $\tprf$, and try to argue that in some models of~$\tprf$, $R^I$ is not valid. But here the
weakness of~$\tprf$ that we were hoping to exploit becomes our worst enemy: $I$ may involve formulas of arbitrary
high quantifier complexity that may potentially denote very complicated combinatorial properties, and we just have no
handle on how to understand them. What we need is that definable sets have manageable structure.

Ideally, we would like to extend $\tprf$ to a (consistent) theory~$T$ with full quantifier elimination. Now, a moment
of reflection tells us that any possible configuration of finitely many functions on a finite set may be realized by
suitable recursive functions, and as such should embed into a model of~$T$. By compactness, \emph{any} $L_\sprf$-structure
should embed in a model of~$T$, thus if such a $T$ exists, it is unique: $T$ must be the model completion of the empty
theory in~$L_\sprf$. (By the empty theory, we mean the theory with no extra-logical axioms.) 

The model completion of the empty $L$-theory is well known and well understood for finite \emph{relational}
languages~$L$. The theory can be axiomatized by a transparent set of ``extension axioms'', and it coincides with
the set of all formulas that hold in random finite $L$-structures with asymptotic probability~$1$. The theory is
$\omega$-categorical, and its unique countable model is the countable random $L$-structure generalizing the Erd\H
os--R\'enyi--Rado random graph; alternatively, it can be described as the Fra\"\i ss\'e limit of the class of all
finite $L$-structures.

It is much less known that the model completion of the empty $L$-theory does, indeed, exist for arbitrary
languages~$L$, as we need here. This was proved by Winkler~\cite{winkler} as a corollary of more general results on
model companions of Skolem expansions of model-complete theories. Note that for languages with functions, most of the
above-mentioned properties of the theory of the random relational structure break down: first-order logic with
functions has no 0--1 law (or even limit law) on finite structures, and there does not seem to be a sensible way of
defining a probability distribution on functions on infinite sets; we will see that the model completion is not
$\omega$-categorical, and not locally finitely satisfiable.

We will now give a self-contained argument that the model completion exists, including an explicit axiomatization by
extension formulas; this will also help us later to determine (syntactically) what open formulas are consistent with the theory. We will denote the model completion as~$\ec_L$, as its models are exactly the (absolutely) existentially
closed $L$-structures.

\begin{Def}\label{def:flas}
Let $L$ be a finite language, and $\Theta$ be a finite set of $L$-terms closed under subterms such that the variables
in~$\Theta$ are among
$x_0,\dots,x_{n-1},y_0,\dots,y_{m-1}$. Let $\ep^R_{\vc t}\in\{0,1\}$ for every $k$-ary relation $R\in L$, and every
$t_0,\dots,t_{k-1}\in\Theta$. Let $\sim$ be an equivalence relation on~$\Theta$ such that:
\begin{enumerate}
\item\label{item:1}
If $R\in L$ is $k$-ary, and $t_i\sim s_i$ for each $i<k$, then $\ep^R_{\vc t}=\ep^R_{\vc s}$.
\item\label{item:2}
If $F\in L$ is $k$-ary, and $t\eqs F(\vc t)\in\Theta$ and $s\eqs F(\vc s)\in\Theta$ satisfy $t_i\sim s_i$ for each $i<k$, then $t\sim s$.
\end{enumerate}
Then the \emph{elementary existential formula} $\exists\vc y\,\theta_{\Theta,\sim,\ep}(\vc x,\vc y)$ is defined by
\begin{equation}\label{eq:1}
\theta_{\Theta,\sim,\ep}(\vc x,\vc y)\eqs\ET_{\substack{t,s\in\Theta\\t\sim s}}t=s\land\ET_{\substack{t,s\in\Theta\\t\nsim s}}t\ne s\land
  \ET_{\substack{R\in L\\\vc t\in\Theta}}R^{\ep^R_{\vc t}}(\vc t),
\end{equation}
where $\psi^1\eqs\psi$, $\psi^0\eqs\neg\psi$.
\end{Def}
\begin{Lem}\label{lem:disj}
Every existential formula $\exists\vc y\,\theta(\vc x,\vc y)$ in a finite language~$L$ is equivalent to a disjunction of elementary
existential formulas with the same free and bound variables.
\end{Lem}
\begin{Pf}
Let $\Theta$ be the set of all subterms of $\theta(\vc x,\vc y)$, and $\Phi$ the set of all (finitely many) atomic
formulas using terms from~$\Theta$. We can write~$\theta$ in full disjunctive normal form in atoms~$\Phi$, and switch
disjunctions with existential quantifiers. Each disjunct has the form~\eqref{eq:1}, except for the conditions on $\sim$
and~$\ep$. However, it is easy to see that if $\sim$ is not an equivalence relation, or if \ref{item:1} or~\ref{item:2}
is violated, then $\theta_{\Theta,\sim,\ep}$ is contradictory.
\end{Pf}
\begin{Def}\label{def:star}
Let $L,\Theta,\vc x,\vc y,\sim,\ep$ be as in Definition~\ref{def:flas}. We define a subset $\Xi\sset\Theta$, and for each $t\in
\Xi$ a term $t^*(\vc x)$, as follows:
\begin{enumerate}
\item\label{item:3}
Every variable $x_i$ is in~$\Xi$, and $x_i^*\eqs x_i$.
\item\label{item:4}
If $t\sim s\in \Xi$, then $t\in \Xi$, and $t^*\eqs s^*$.
\item\label{item:5}
If $t\eqs F(t_0,\dots,t_{k-1})\in\Theta$, and $t_0,\dots,t_{k-1}\in \Xi$, then $t\in \Xi$, and $t^*\eqs F(t_0^*,\dots,t_{k-1}^*)$.
\end{enumerate}
If more than one clause applies to put $t\in \Xi$, we define $t^*$ using any of them; the choice does not matter. The
only relevant property is the following observation:
\begin{equation}\label{eq:2}
\vdash\theta_{\Theta,\sim,\ep}(\vc x,\vc y)\to t(\vc x,\vc y)=t^*(\vc x)
\end{equation}
for every $t\in \Xi$. Finally, we define an open formula $\theta_{\Theta,\sim,\ep}^*(\vc x)$ as
\[\ET_{\substack{t,s\in\Xi\\t\sim s}}t^*=s^*
   \land\ET_{\substack{t,s\in \Xi\\t\nsim s}}t^*\ne s^*\land
  \ET_{\substack{R\in L\\\vc t\in\Xi}}R^{\ep^R_{\vc t}}(\vc t^*)\land
  \ET_{\substack{t\eqs F(\vc t)\in\Xi\\\vc t\in\Xi}}t^*=F(\vc t^*).\]
\end{Def}
\begin{Prop}\label{prop:res}
Let $L,\Theta,\sim,\ep$ be as in Definition~\ref{def:flas}, $M$ an $L$-structure, and $\vc u\in M$. The following are
equivalent.
\begin{enumerate}
\item\label{item:6}
There exists an extension $N\Sset M$ such that $N\model\exists\vc y\,\theta_{\Theta,\sim,\ep}(\vc u,\vc y)$.
\item\label{item:7}
$M\model\theta_{\Theta,\sim,\ep}^*(\vc u)$.
\end{enumerate}
\end{Prop}
\begin{Pf}
\ref{item:6}${}\to{}$\ref{item:7} follows immediately from the definitions and property~\eqref{eq:2}.

\ref{item:7}${}\to{}$\ref{item:6}: Put $N=M\cupd(\Theta\bez\Xi)/{\sim}$, and identify $t/{\sim}$ with $t^*(\vc u)$ for
$t\in\Xi$. If $t_0,\dots,t_{k-1}\in\Theta$, $R\in L$, and $t=F(\vc t)\in\Theta$, define
\begin{align*}
R^N(t_0/{\sim},\dots,t_{k-1}/{\sim})&\iff\ep^R_{\vc t}=1,\\
F^N(t_0/{\sim},\dots,t_{k-1}/{\sim})&=t/{\sim}.
\end{align*}
Using $M\model\theta_{\Theta,\sim,\ep}^*(\vc u)$ and the properties of $\sim$ and~$\Xi$, it is easy to check that the
definition is independent of the choice of representatives, and agrees with~$M$ if $t_0,\dots,t_{k-1}\in\Xi$. 
We can thus extend the definition with the original structure of~$M$, and then arbitrarily to the remaining tuples
from~$N$. The definition ensures
\[N\model\theta_{\Theta,\sim,\ep}(\vc x/{\sim},\vc y/{\sim}),\]
where $x_i/{\sim}=u_i$.
\end{Pf}

In other words, Proposition~\ref{prop:res} shows that $\theta_{\Theta,\sim,\ep}^*$ is equivalent to the \emph{resultant} (see
\cite[\S7.2]{hodges:sh}) of the elementary existential formula $\exists\vc y\,\theta_{\Theta,\sim,\ep}(\vc x,\vc y)$
in the empty theory.

\begin{Lem}\label{lem:finite}
Let $M$ be an $L$-structure.
\begin{enumerate}
\item\label{item:8}
If $M$ is existentially closed, then so is $M\res L'$ for every $L'\sset L$.
\item\label{item:14}
If $M\res L'$ is e.c.\ for every finite $L'\sset L$, then $M$ is e.c.
\end{enumerate}
\end{Lem}
\begin{Pf}
\ref{item:8}: Every extension of $M\res L'$ can be expanded to an extension of~$M$. \ref{item:14}: Assume that $M\sset
N\model\fii(\vc u)$, where $\vc u\in M$, and $\fii$ is an existential formula. Let $L'\sset L$ be a finite
sublanguage containing all symbols occurring in~$\fii$. Then $\fii(\vc u)$ holds in $N\res L'$, hence in $M\res L'$
(and~$M$) as $M\res L'$ is e.c.
\end{Pf}
\begin{Def}\label{def:ec}
If $L$ is a finite language, let $\ec_L$ denote the theory axiomatized by the formulas
\[\theta_{\Theta,\sim,\ep}^*(\vc x)\to\exists y\,\theta_{\Theta,\sim,\ep}(\vc x,y)\]
for all $\Theta,\sim,\ep$ as in Definition~\ref{def:flas} with $m=1$.

For infinite~$L$, we put $\ec_L=\bigcup\{\ec_{L'}:L'\sset L\text{ finite}\}$.
\end{Def}
\pagebreak[2]
\begin{Thm}\label{thm:comp}
For any language~$L$, $\ec_L$ is a model completion of the empty $L$-theory:
\begin{enumerate}
\item\label{item:9}
$\ec_L$ has elimination of quantifiers.
\item\label{item:10}
Models of $\ec_L$ are exactly the existentially closed $L$-structures; in particular, every $L$-structure embeds in a
model of~$\ec_L$. 
\end{enumerate}
Moreover, if $L$ is recursively presented, then $\ec_L$ is decidable.
\end{Thm}
\begin{Pf}
By Lemma~\ref{lem:finite}, we may assume that $L$ is finite.

\ref{item:9}: Proposition~\ref{prop:res} implies the converse implications
\[{}\vdash\exists y\,\theta_{\Theta,\sim,\ep}(\vc x,y)\to\theta_{\Theta,\sim,\ep}^*(\vc x)\]
for each axiom of~$\ec_L$. In view of Lemma~\ref{lem:disj}, this shows that every formula with one existential
quantifier is equivalent to an open formula over~$\ec_L$, hence the same follows for all formulas by induction on
complexity.

\ref{item:10}: Proposition~\ref{prop:res} implies that any e.c.\ model validates~$\ec_L$. The converse would also follow from
Proposition~\ref{prop:res}, were it not for the restriction to $m=1$ in the definition of~$\ec_L$. However, if $\fii(\vc x)$ is
an existential formula, $N\Sset M\model\ec_L$, and $N\model\fii(\vc u)$ for some $\vc u\in M$, we can assume $N$ is
e.c.\ by extending it further if necessary. Thus, $N\model\ec_L$, and $M$ is an elementary substructure of~$N$
by~\ref{item:9}, hence $M\model\fii(\vc u)$.

As for decidability, $\ec_L$ is clearly r.e., hence quantifier elimination is effective (in fact, the argument above
gives an explicit algorithm). Now, if $\fii$ is a quantifier-free (or even universal) sentence, then $\ec_L\vdash\fii$
iff $\nul\vdash\fii$ iff $\fii$ holds in all finite models of cardinality bounded by the number of subterms $t$
of~$\fii$: if $M\nmodel\fii$, let $M_0=\{t^M:t\text{ a subterm of }\fii\}\sset M$. We define realizations of relation
and function symbols from~$L$ in~$M_0$ to be the same as in~$M$, except when a function symbol takes a value
outside~$M_0$, in which case we redefine it as an arbitrary element of~$M_0$. Then $M_0\nmodel\fii$, since
$t^{M_0}=t^M$ for all subterms of~$\fii$.
\end{Pf}

\begin{Cor}\label{cor:dec}
If $M$ is a recursive $L$-structure, $\ec_L+\diag(M)$ is a decidable complete theory.
\noproof\end{Cor}

\section{Interpretability in existential theories}\label{sec:interpr-exist-theor}

As explained in the beginning of Section~\ref{sec:model-compl-empty}, our intention for discussing $\ec_L$
is that we want to show noninterpretability of certain theories~$S$ in $\tprf$ (which is an existential---in fact,
quantifier-free---theory) by showing their noniterpretability in completions of $\ec_L$. Now, it is not a priori clear if
this strategy is any good: why should we expect that $S$ is, indeed, not intepretable in an extension of~$\ec_L$? After
all, $\ec_L$ is a fairly nontrivial theory, hence it may interpret~$S$ even if no existential theory can; to begin
with, it interprets~$\ec_L$ itself.

To allay our fears, we will prove in the present section that a large class of theories~$S$ is immune to such
shenanigans: specifically, a theory axiomatized by \emph{$\exists\forall$~sentences} is interpretable in a completion
of some $\ec_L$ \emph{if and only if} it is interpretable in a consistent quantifier-free or existential theory. We
believe this characterization to be of independent interest, which is why we prove it in detail even though it is only of
indirect relevance for our original goal, in that it assures us that we set off in the right direction; of
course, we would eventually find that anyway when we got to the point.

We start with a few auxiliary lemmas whose basic intention is that if a theory is (weakly) interpretable in $\ec_L$, we
can make $L$ and the interpretation ``nice''.
\pagebreak[2]
\begin{Lem}\label{lem:intp}
\ \begin{enumerate}
\item\label{item:11}
If $L'\Sset L$, $\ec_{L'}$ is a conservative extension of~$\ec_L$.
\item\label{item:34}
If $L'=L\cup\{c_i:i\in I\}$, where $c_i\notin L$ are constants, then $\ec_{L'}$ is axiomatized by $\ec_L$. In
particular, $\ec_{L'}\vdash\fii(\vc c)$ iff $\ec_L\vdash\forall\vc x\,\fii(\vc x)$.
\item\label{item:12}
If $L$ contains a constant~$c$, and $F,R\notin L$ are a $k$-ary function and relation symbol (resp.), then
\[R^I(\vc x)\eqs F(\vc x)=c\]
provides a faithful interpretation of $\ec_{L\cup\{R\}}$ in~$\ec_{L\cup\{F\}}$.
\item\label{item:13}
If $L$ contains an at least binary function, or at least two unary functions, then $\ec_L$ has pairing.
\end{enumerate}
\end{Lem}
\begin{Pf}
\ref{item:11}: $\ec_{L'}\Sset\ec_L$ follows from Lemma~\ref{lem:finite}. On the other
hand, any model $M\model\ec_L$ has an expansion to an $L'$-structure $M'$, which has an extension $N'\model\ec_{L'}$.
$M$ is an elementary substructure of~$N'\res L$ by quantifier elimination for~$\ec_L$.

\ref{item:34} follows from~\ref{item:11} and the fact that an expansion of an e.c.\ model by constants is e.c.

\ref{item:12}: Let $M\model\ec_{L\cup\{F\}}$, and $\fii(\vc u)$ be an existential $(L\cup\{R\})$-formula with $\vc
u\in M$, satisfied in an extension $N\Sset M^I$. We can extend $F^M$ to $N$ so that $F^N(\vc a)=c$ iff $R^N(\vc a)$
for all $\vc a\in N$. Then $\p{N\res L,F^N}\model\fii^I(\vc u)$, where $\fii^I$ is an existential formula, hence
$M\model\fii^I(\vc u)$ by e.c., and $M^I\model\fii(\vc u)$. Thus, $M^I\model\ec_{L\cup\{R\}}$, which shows that
$\ec_{L\cup\{F\}}\vdash\ec_{L\cup\{R\}}^I$. On the other hand, assume $\ec_{L\cup\{F\}}\vdash\fii^I$, and
$M\model\ec_{L\cup\{R\}}$. Let $F^M\colon M^k\to M$ be such that $F^M(\vc a)=c^M$ iff $R^M(\vc a)$, and
$N\Sset\p{M\res L,F^M}$ be a model of~$\ec_{L\cup\{F\}}$. We have $N^I\model\ec_{L\cup\{R\}}+\fii$, and $M\preceq N^I$
by quantifier elimination, hence $M\model\fii$.

\ref{item:13}: If $L$ contains two distinct unary function symbols $L(x)$, $R(x)$, we claim that
\[\ec_L\vdash\forall x,y\,\exists z\,(L(z)=x\land R(z)=y).\]
Let $M\model\ec_L$, and $a,b\in M$. Put $N=M\cupd\{c\}$, where $L^N(c)=a$, $R^N(c)=b$, and the realization of other
functions or relations on tuples involving~$c$ is arbitrary. Then $N\model\exists z\,(L(z)=a\land R(z)=b)$, hence the
same holds in~$M$ by e.c.

If $L$ contains a $k$-ary function $F$ for $k\ge2$, we may use a similar argument with e.g.
\begin{align*}
L(x)&=F(F(x,\dots,x),x,\dots,x),\\
R(x)&=F(x,\dots,x,F(x,\dots,x)),
\end{align*}
and $N=M\cupd\{c,F(c,\dots,c)\}$.
\end{Pf}

\begin{Def}\label{def:int-qf}
We will call a parameter-free translation \emph{quantifier-free} if its domain, and the translations of all predicate
symbols as well as equality are given by quantifier-free formulas, and the translations of all function symbols are
given piecewise by terms, where the pieces are finitely many and quantifier-free definable.
\end{Def}

Recall Lemma~\ref{lem:coi-wi}.

\begin{Lem}\label{lem:binary}
For any language~$L$, there is a quantifier-free one-piece one-dimensional parameter-free unrelativized cointerpretation with
absolute equality of $\ec_{L_2}$ in~$\ec_L$, where $L_2$ consists of a single binary function and at most~$\lh L$
constants. If $L$ is countable, one constant suffices.
\end{Lem}
\begin{Pf}
Using Lemma~\ref{lem:intp}~\ref{item:12} (which may be applied in parallel to all relations using the same argument), we may assume $L$ contains no relations. Let $L_2$ be the language consisting of
a binary function $(x,y)$, the constants of~$L$, and new constants $c_F$ for every nonconstant function $F\in L$. For
$n\ge1$, write
\[(x_0,\dots,x_{n-1}):=(x_0,(x_1,\dots,(x_{n-2},x_{n-1})\dots)).\]
Let $I$ be the translation of $L$ into~$L_2$ defined by $c^I=c$ for constants $c\in L$, and
\[F^I(x_0,\dots,x_{n-1})=((c_F,x_0,\dots,x_{n-1}),x_0,\dots,x_{n-1},x_0,\dots,x_{n-1})\]
for $n$-ary functions $F\in L$, $n>0$. Let
\[\ec_{L_2}^*=\ec_{L_2}+\{c_F\ne c_G:F\ne G\in L\}.\]
We will show that $I$ is a faithful interpretation of $\ec_L$ in~$\ec_{L_2}^*$, which implies it is also a
cointerpretation of $\ec_{L_2}$ in~$\ec_L$.
\begin{Cl}\label{cl:extend}
\ \begin{enumerate}
\item\label{item:36}
If $N$ is an $L_2$-structure such that the constants $c_F^N$ are pairwise distinct, and $M$ is an extension of the
$L$-structure $N^I$, there is an extension $K\Sset N$ such that $K^I\Sset M$.
\item\label{item:35}
If $M$ is an $L$-structure, there is an $L_2$-structure $N$ such that $N^I\Sset M$, and the constants $c_F^N$ are
distinct.
\end{enumerate}
\end{Cl}
\begin{Pf*}
\ref{item:36}: Let $K$ be the disjoint union $M\cupd M^{<\omega}$, with constants realized as in~$N$, and
\begin{align*}
(a,b)^K&=(a,b)^N&&a,b\in N,\\
(a_0,\p{a_1,\dots,a_n})^K&=\p{a_0,\dots,a_n}&&a_0\in M,\\
(\p{c_F,u_0,\dots,u_{i-1}},\p{a_0,\dots,a_{n-1},v_0,\dots,v_{j-1}})^K&=F^M(a_0,\dots,a_{n-1})&&\text{$F$ $n$-ary, $i,j\ge0$,}\\
(a,b)^K&=\p{a}&&a,b\in M,\{a,b\}\nsset N,\\
(a,b)^K&=\p{}&&\text{all other cases.}
\end{align*}
We need to check that if $F\in L$ is $n$-ary with $n>0$, and $a_0,\dots,a_{n-1}\in M$, then
\[F^M(\vc a)=(F^I(\vc a))^K.\]
We may assume $\vc a\notin N^n$, as otherwise the statement follows from $N^I\sset M$. Let $i<n$ be maximal such
that $a_i\notin N$, and $i'=\min\{i,n-2\}$. It follows from the definition that
\begin{align*}
(a_0,\dots,a_{n-1},a_0,\dots,a_{n-1})^K&=\p{a_0,\dots,a_{n-1},a_0,\dots,a_{i'}},\\
(c_F,a_0,\dots,a_{n-1})^K&=\p{c_F,a_0,\dots,a_{i'}},
\end{align*}
hence
\[((c_F,a_0,\dots,a_{n-1}),a_0,\dots,a_{n-1},a_0,\dots,a_{n-1})^K=F^M(a_0,\dots,a_{n-1})\]
as required.

\ref{item:35}: By extending $M$ if necessary, we may assume $\lh M\ge\lh L$, hence we can fix pairwise distinct
elements $c_F^N\in M$. Put $N=M\cupd M^{<\omega}$, and define
\begin{align*}
(a,b)^N&=\p{a,b}&&a,b\in M,\\
(a_0,\p{a_1,\dots,a_n})^N&=\p{a_0,\dots,a_n}&&a_0\in M,\\
(\p{c_F,u_0,\dots,u_{i-1}},\p{a_0,\dots,a_{n-1},v_0,\dots,v_{j-1}})^K&=F^M(a_0,\dots,a_{n-1})&&\text{$F$ $n$-ary, $i,j\ge0$,}\\
(a,b)^N&=\p{}&&\text{otherwise.}
\end{align*}
We have $M\sset N^I$ by a similar (but easier) argument as in~\ref{item:36}.
\end{Pf*}

In order that $I$ interprets $\ec_L$ in~$\ec_{L_2}^*$, it suffices to show that if
$N\model\ec_{L_2}^*$, then $N^I$ is e.c. Now, if an existential $L$-formula $\fii$ with parameters from~$N$ is
satisfiable in $M\Sset N^I$, then $M\sset K^I$ for some $K\Sset N$ by the claim, which thus satisfies the existential formula
$\fii^I$. It follows that $N\model\fii^I$ as $N$ is e.c., i.e., $N^I\model\fii$.

To show that $I$ is faithful, let $\ec_{L_2}^*\vdash\fii^I$, and $M\model\ec_L$. By the claim, there is $N$ with
the elements $c_F^N$ pairwise distinct such that $N^I\Sset M$. By extending it if necessary, we may assume $N\model\ec_{L_2}^*$,
hence $N^I\model\fii$. Also, $N^I\model\ec_L$, hence $M\preceq N^I$ by quantifier elimination, which gives
$M\model\fii$.

Finally, let $L$ be countable, and enumerate it as $\{F_k:k\in\omega\}$. Let $L_2$ be the language consisting of
$(x,y)$ and a single constant~$c$. We modify the construction above as follows: we employ the closed terms
\[c_k=(\dots((\underbrace{c,c),c),\dots,c}_{\text{$n+2$ times}}).\]
in place of $c_{F_k}$, and if $F_k$ is a constant, we put
$F_k^I=(c_k,c_k)$. (In particular, we redefine $\ec_{L_2}^*$ to state that all the $c_k$ are pairwise
distinct.) Then it is easy to check that the argument still goes through: the only place where the exact composition
of~$c_F$ matters is in the proof of part~\ref{item:35} of the claim, and we can fix it e.g.\ by making $N$ the set of
all finite binary trees with leaves labelled by $M\cupd\{c\}$, where $a\in M$ is identified with a one-node tree,
and $(x,y)^N$ is the tree whose root has children $x,y$, except for
\begin{align*}
((c_k,a_0,\dots,a_{n-1})^N,(a_0,\dots,a_{n-1},a_0,\dots,a_{n-1})^N)^N&=F_k^M(a_0,\dots,a_{n-1})&&\text{$F_k$ $n$-ary, $\vc a\in M$,}\\
(c_k^N,c_k^N)^N&=F_k^M&&\text{$F_k$ constant,}
\end{align*}
as needed to make the interpretation work.
\end{Pf}

\begin{Cor}\label{cor:simple}
If a theory~$T$ is weakly interpretable in~$\ec_L$ for some~$L$, it has a one-piece one-dimensional parameter-free unrelativized
interpretation in a consistent extension of some $\ec_{L_2}$, where $L_2$ consists of a binary function, and at most
$\lh L$ constants. If $L$ is countable, one constant suffices.
\end{Cor}
\begin{Pf}
We can make the interpretation one-piece as $\ec_L$ proves there are at least two elements.
We can assume $L$ contains a constant by Lemma~\ref{lem:intp}~\ref{item:11}, and that it is purely functional by~\ref{item:12}. Then we can make the interpretation
one-dimensional by~\ref{item:13}, and parameter-free by expanding $L$ with constants for the parameters,
using~\ref{item:34}. We can also assume to
have a constant~$c$ denoting an element in the domain of the interpretation, and then it is easy to make
the interpretation unrelativized by equating (i.e., extending the interpreted equality) elements outside the domain
with~$c$. Finally, we can compose the interpretation with the one from Lemma~\ref{lem:binary} to make the language as
needed.
\end{Pf}

Note that the argument in Corollary~\ref{cor:simple} does not guarantee that the interpretation is quantifier-free: while the
domain and the translations of all symbols can be made quantifier-free formulas just by quantifier elimination, this
does not ensure function symbols are given piecewise by terms. This will in fact pose a serious challenge in the proof
of the characterization below, and we will need results on elimination of imaginaries from Appendix~\ref{sec:elim-imag} to
deal with it.

\begin{Thm}\label{thm:ea}
Let $T$ be an $\exists\forall$-axiomatized theory in a language~$L_T$. The following are
equivalent.
\begin{enumerate}
\item\label{item:16}
$T$ is interpretable in a consistent existential theory.
\item\label{item:17}
$T$ has a quantifier-free interpretation~$I$ in a consistent quantifier-free theory~$S$ such that $I$ and the language
of~$S$ obey the conditions in Corollary~\ref{cor:simple}, except that $I$ may be multi-dimensional if $L_T$ contains a proper
function symbol.
\item\label{item:18}
$T$ is weakly interpretable in $\ec_L$ for some language~$L$, w.l.o.g.\ obeying the same conditions as in~\ref{item:17}.
\end{enumerate}
If $L_T$ is finite, and $T$ is recursively axiomatized, we can also make the interpreting theories
recursively axiomatized.
\end{Thm}
\begin{Pf}
\ref{item:17}${}\to{}$\ref{item:16} is trivial, and \ref{item:16}${}\to{}$\ref{item:18} follows from the fact that
every consistent existential theory is consistent with~$\ec_L$ in the same language by Theorem~\ref{thm:comp}.

\ref{item:18}${}\to{}$\ref{item:17}: By expanding $L$ and~$L_T$ with Henkin constants for the existential quantifiers
in axioms of~$T$ using Lemma~\ref{lem:intp}~\ref{item:34}, we may assume that $T$ is universal.
By Corollary~\ref{cor:simple}, $T$ has a one-piece one-dimensional parameter-free unrelativized interpretation~$J$ in a consistent
theory $\ec_L+S$, where $L$ consists of a binary function and constants. By quantifier elimination, we may assume $S$
is a set of quantifier-free sentences, and the $J$-translations of equality and all symbols of~$L_T$ are given
by quantifier-free formulas. By expanding the language~$L$ further, we may assume
that constants of~$L_T$ are interpreted by constants (or constant terms) of~$L$. In the countable case,
we may apply Lemma~\ref{lem:binary} again to reduce the number of constants to one.

If $L_T$ contains proper function symbols, we need more work, as we cannot add Skolem functions in the same way as
constants. As we will explain in Appendix~\ref{sec:elim-imag}, $\ec_L$ has weak elimination of imaginaries, and as a corollary, we obtain in
Proposition~\ref{prop:equi-expl} an explicit description of definable equivalence relations that we apply to $=^J$.
Using~\eqref{eq:5}, we see that the collection of equivalence classes that make up the domain of~$J$ can be definably
split in finitely many pieces, where the $i$-th piece is in definable bijection with a collection of $m_i$-element sets
of $r$-tuples (represented by an equivalence relation on $m_ir$-tuples as in Definition~\ref{def:ei}). The upshot is that we may replace
$J$ with an equivalent piece-wise interpretation~$I$ that almost has absolute equality, in the sense that all
equivalence classes of $=^I$ have bounded finite size. Consequently, the translation $F^I$ of any function symbol $F\in
L_T$, when viewed as a relation on tuples rather than on their equivalence classes, is a total multifunction with only
finitely many values. By Lemma~\ref{lem:acl} and a compactness argument, there is a piecewise term-definable function that
picks one possible value of such a multivalued function. Thus, $I$ is a quantifier-free interpretation. Since we may
assume $L$ includes a pair of constants $c,d$ such that $S\vdash c\ne d$, we can make $I$ a one-piece interpretation;
it is still parameter-free, and we can make it unrelativized as above, but it may be multi-dimensional.
(Lemma~\ref{lem:intp} does not give a pairing \emph{function}, hence it is unclear if we can make the interpretation
one-dimensional without sacrificing the property that translations of functions are piecewise term-definable.)

The result of these manipulations is that $T^I$ is a \emph{universal} subtheory of $\ec_L+S$, as we made sure all
existential quantifiers needed are witnessed (piecewise) by terms. Thus, $T^I$ is in fact included in~$S$, i.e., $I$ is
an interpretation of $T$ in~$S$ which satisfies all the requirements.

Finally, let $T$ be an r.e.\ theory in a finite language. We have shown that if $T$ is interpretable in a consistent
existential theory, there is a quantifier-free unrelativized one-piece parameter-free interpretation of~$T$ in a consistent
extension of~$\ec_{L_2}$, where $L_2$ consists of a constant and a binary function. (The interpretation is
automatically recursive, as the language is finite.) The universal Henkin expansion~$T^H$
of~$T$ is still r.e., and we can assign the Henkin constants in a recursive way to new constants added to~$L_2$ so that
we get an interpretation~$I$ with the same properties of~$T^H$ in a consistent extension of~$\ec_L$, where $L$
consists of a binary function and countably many constants, and $I$ is recursive. The cointerpretation from
Lemma~\ref{lem:binary} is also recursive, hence we can reduce the language back to~$L_2$. Then $\ec_{L_2}+(T^H)^I$ is an
r.e.\ theory, hence by effectiveness of quantifier elimination, it is equivalent to $\ec_{L_2}+S$ for an r.e.\
quantifier-free $L_2$-theory~$S$. By the argument above, $I$ is an interpretation of $T$ in~$S$, as $(T^H)^I$ is a
universal theory.
\end{Pf}
\pagebreak[2]
\begin{Rem}
Theorem~\ref{thm:ea} does not extend to $\forall\exists$ theories~$T$. On the one hand, any theory interpretable in a
consistent
existential theory is locally finitely satisfiable (notice also that any consistent $\exists\forall$ theory in a
relational language is locally finitely satisfiable). On the other hand, $\ec_L$ itself is a $\forall\exists$ theory interpretable
in~$\ec_L$, and if $L$ contains a nonconstant function symbol, then $\ec_L$ is
not locally finitely satisfiable: for example, if we
have a unary function $F(x)$, then $\ec_L$ proves the formula
\[\forall x,y\,\exists z\,(z\ne x\land F(z)=y)\]
with no finite model. 

We note that if $L$ contains only at most unary relations and constants, then $\ec_L$ and any its consistent extension
is an existential theory, and easily seen to be interpretable in~$\tprfu$ for $L$~finite. If $L$ consists of relations
and constants, but is not unary, then $\ec_L$ (i.e., essentially the theory of the random structure) is genuinely
$\forall\exists$, but still locally finitely satisfiable, hence interpretable in~$R$ for $L$ finite by Visser's
Theorem~\ref{thm:r-locfinsat}.
\end{Rem}

\begin{Que}\label{que:reprf}
Is every consistent r.e.\ existential theory interpretable in~$\tprf$?
\end{Que}
\begin{Que}\label{que:erandgrf}
Is the theory of the random graph interpretable in a consistent existential theory?
\end{Que}

\section{Classification of $\ec_L$}\label{sec:classification-ec_l}

We now proceed to the main results of the paper, showing that certain theories are not interpretable in any
existentially axiomatized theory by way of establishing tameness properties of~$\ec_L$. We will mostly deduce
them from the following statement, showing the impossibility of certain configurations in models of~$\ec_L$.

In order to keep the proof self-contained and accessible to wider audience, we will not use any results on
indiscernibles (though they are lurking in our application of Ramsey's theorem).

Recall that a relation $R\sset X^2$ is \emph{asymmetric} if there are no $a,b\in X$ such that $R(a,b)\land R(b,a)$.
\begin{Thm}\label{thm:main}
For any language~$L$ and formula $\fii(\vc z,\vc x,\vc y)$ with $\Lh(\vc x)=\Lh(\vc y)$, there is a constant~$n$
with the following property. Let $M\model\ec_L$ and $\vc u\in M$ be such that
\[M\model\exists\vc x_0,\dots,\vc x_{n-1}\,\ET_{i<j<n}\fii(\vc u,\vc x_i,\vc x_j).\]
Then for every~$m\in\omega$ and every asymmetric relation~$R$ on $\{0,\dots,m-1\}$,
\[M\model\exists\vc x_0,\dots,\vc x_{m-1}\,\ET_{\p{\alpha,\beta}\in R}\fii(\vc u,\vc x_\alpha,\vc x_\beta).\]
\end{Thm}
\begin{Pf}
By Theorem~\ref{thm:comp} and Lemma~\ref{lem:intp}, we may assume $L$ contains no relations, all the tuples have length one, and $\fii$ is open. Let $\tau$ be the
number of subterms of~$\fii$, and $\tau^*=2^{256\tau^2}$. Using Ramsey's theorem, let $n$ be sufficiently large so that
\[n\to(7)^4_{\tau^*}.\]
Fix $M\model\ec_L$, $u\in M$, and $\{a_i:i<n\}\sset M$ such that $M\model\fii(u,a_i,a_j)$ for $i<j<n$. In order to
simplify the notation, we will assume $u$ is given by a constant of~$L$, and write just $\fii(x,y)$; this does not
increase the number of subterms of~$\fii$. Let $S$ be the set of all subterms $t(x,y)$ of~$\fii$, and for every
$i_0<i_1<i_2<i_3<n$, define
\[\tp(i_0,i_1,i_2,i_3)=\{\p{u_0,u_1,u_2,u_3,t,s}\in\jvn4^4\times S^2:
  t^M(a_{i_{u_0}},a_{i_{u_1}})=s^M(a_{i_{u_2}},a_{i_{u_3}})\}.\]
Since $\lh S\le\tau$, $\tp$ is a colouring of quadruples of numbers below~$n$ by at most $\tau^*$
colours. Thus, we can find a $7$-element homogeneous set $H\sset\{0,\dots,n-1\}$ for~$\tp$; without loss of
generality $H=\{0,\dots,6\}$.

Fix a set of variables $\{y_\alpha:\alpha<m\}$, and put
\[\Theta=\{t(y_\alpha,y_\beta):\p{\alpha,\beta}\in R,t(x,y)\in S\}.\]
If $t\in\Theta$, let $V(t)$ denote the set of $\alpha<m$ such that $y_\alpha$ occurs in~$t$; note that $\lh{V(t)}\le2$.
A \emph{realization} of~$t$ is an injective mapping $r\colon V(t)\to H$ such that
\[\alpha,\beta\in V(t)\land\p{\alpha,\beta}\in R\implies r(\alpha)<r(\beta).\]
Notice that this condition is void if $t$ depends on at most one variable; otherwise it concerns a unique pair
$\p{\alpha,\beta}$. If $r$ is a realization of~$t$, let $r(t)\in M$ be the value of the term
resulting from~$t$ by replacing each variable $y_\alpha$ with $a_{r(\alpha)}$.

A \emph{joint realization} of a set of terms $\{t_0,\dots,t_{k-1}\}\sset\Theta$ is an injective mapping $r\colon V(t_0)\cup\dots\cup V(t_{k-1})\to H$ such that
$r\res V(t_i)$ is a realization of~$t_i$ for $i<k$. Note that any pair $\{t,s\}\sset\Theta$ has a joint
realization, as $R$ has no cycles of length at most~$2$.

If $t,s\in\Theta$, and $r$ is a joint realization of $t$ and~$s$, we define
\[t\sim s\iff r(t)=r(s).\]
\begin{Cl}\label{cl:welldef}
The definition of~$\sim$ is independent of the choice of~$r$.
\end{Cl}
\begin{Pf*}
First, if two joint realizations $r,r'$ satisfy
\begin{equation}\label{eq:16}
r(\alpha)<r(\beta)\iff r'(\alpha)<r'(\beta)
\end{equation}
for all $\alpha,\beta\in V(t)\cup V(s)$, then
\begin{equation}\label{eq:10}
r(t)=r(s)\iff r'(t)=r'(s)
\end{equation}
by homogeneity for~$\tp$. This condition holds automatically if
\begin{itemize}
\item $V(t)\sset V(s)$ or $V(s)\sset V(t)$, or
\item $V(t)=\{\alpha,\beta\}$, $V(s)=\{\beta,\gamma\}$, where $\p{\alpha,\beta},\p{\beta,\gamma}\in R$, or vice versa.
\end{itemize}
Assume $t=t(y_\alpha,y_\beta)$, $s=s(y_\alpha,y_\gamma)$, where $\beta\ne\gamma$, and
$\p{\alpha,\beta},\p{\alpha,\gamma}\in R$ (the case with $\p{\beta,\alpha},\p{\gamma,\alpha}\in R$ is symmetric). By~\eqref{eq:10}, it suffices
to consider the case where $r(\alpha)=r'(\alpha)=0$, $r(\beta)=1$, $r(\gamma)=r'(\gamma)=2$, $r'(\beta)=3$.
Using~\eqref{eq:10}, we have
\[t^M(a_0,a_1)=s^M(a_0,a_2)\implies t^M(a_0,a_3)=s^M(a_0,a_4)=t^M(a_0,a_1)=s^M(a_0,a_2),\]
and the converse implication is symmetric.

The remaining case is when $V(t)$ and $V(s)$ are disjoint and nonempty. It suffices to show that if $r(t)=r(s)$ for
some joint realization~$r$, there is a constant $a\in M$ such that $r'(t)=a$ for every realization $r'$ of~$t$ (whence
the same holds for~$s$ by symmetry). Assume $t$ depends on two variables $y_\alpha,y_\beta$ with $\p{\alpha,\beta}\in R$ (the
unary case is easier). Using~\eqref{eq:10}, we may assume that the realization $r_1(\alpha)=r(\alpha)-1$,
$r_1(\beta)=r(\beta)+1$ of $t$ is within bounds, and disjoint from $r(V(s))$. Then $r'=r_1\cup(r\res V(s))$ is a joint
realization of $\{t,s\}$ such that \eqref{eq:16} holds, hence using \eqref{eq:10} again, it follows that
\[t^M(a_{r(\alpha)-1},a_{r(\beta)+1})=r(s)=t^M(a_{r(\alpha)},a_{r(\beta)}).\]
Applying homogeneity, we have
\[t^M(a_i,a_j)=t^M(a_0,a_6)=t^M(a_k,a_l)\]
for every $0<i<j<6$, $0<k<l<6$. Since every set $\{i,j,k,l\}$ is order-isomorphic to some not involving $0,6$, we
obtain
\[t^M(a_i,a_j)=t^M(a_k,a_l)\]
for all $i<j$, $k<l$ using homogeneity again.
\end{Pf*}

Thus, $\sim$ is a well-defined relation on~$\Theta$. It is clearly reflexive and symmetric. If
$t\eqs F(t_0,\dots,t_{k-1})$ and $s\eqs F(s_0,\dots,s_{k-1})$ are in~$\Theta$, and $r$ is a joint realization of $t$ and $s$,
it is also a joint realization of each $\{t_i,s_i\}$, hence
\[t_0\sim s_0,\dots,t_{k-1}\sim s_{k-1}\implies t\sim s.\]
\pagebreak[2]
\begin{Cl}\label{cl:trans}
$\sim$ is transitive.
\end{Cl}
\begin{Pf*}
Assume that $t\sim s\sim u$. If there exists a joint realization $r$ of $\{t,s,u\}$, we immediately obtain
$r(t)=r(s)=r(u)$, hence $t\sim u$. If not, we must have $t=t(y_\alpha,y_\beta)$, $s=s(y_\beta,y_\gamma)$,
$u=u(y_\gamma,y_\alpha)$, where $\p{\alpha,\beta},\p{\beta,\gamma},\p{\gamma,\alpha}\in R$. Applying alternately
$t\sim s$ and $s\sim u$, we obtain
\[t^M(a_3,a_4)=s^M(a_4,a_5)=u^M(a_5,a_6)=s^M(a_1,a_5)=t^M(a_0,a_1)=s^M(a_1,a_2)=u^M(a_2,a_3),\]
hence $r(t)=r(u)$ under the joint realization of $t,u$ such that $r(\gamma)=2$, $r(\alpha)=3$, and $r(\beta)=4$. (This
argument in fact shows that with such a cyclic dependency, the values of all three terms are independent of
the realization.)
\end{Pf*}

Let $\Xi\sset\Theta$ and $\{t^*:t\in\Xi\}$ be as in Definition~\ref{def:star}, for empty~$\vc x$, and (since $L$ has no
relation symbols) empty~$\ep$. By induction on the
definition of~$t\in\Xi$, we see that the value of the closed term $t^*$ in~$M$ coincides with~$r(t)$ for any
realization $r$ of~$t$. This and the definition of~$\sim$ implies that
\[M\model\theta_{\Theta,\sim,\ep}^*,\]
hence by Lemma~\ref{lem:finite} and existential closedness of~$M$,
\[M\model\exists y_0,\dots,y_{m-1}\,\theta_{\Theta,\sim,\ep}(\vc y).\]
If $b_0,\dots,b_{m-1}\in M$ witness this, and $\p{\alpha,\beta}\in R$ and $i<j\in H$, we have
\[M\model\psi(b_\alpha,b_\beta)\iff M\model\psi(a_i,a_j)\]
for every subformula $\psi$ of~$\fii$. It follows that
\[M\model\ET_{\p{\alpha,\beta}\in R}\fii(b_\alpha,b_\beta)\]
as required.
\end{Pf}

We draw two principal conclusions from Theorem~\ref{thm:main}. For the first one, notice that the theory below is
interpretable in the theory~$R$ just by taking $<$ for~$\in$: then \eqref{eq:3} is witnessed by $x_i=\num i $, $z=\num
n$ due to axiom~\eqref{eq:12}.
\begin{Cor}\label{thm:existsinfty}
The theory in the language $\p\in$ axiomatized by the sentences
\begin{equation}\label{eq:3}
\exists z,x_0,\dots,x_{n-1}\,\Bigl(\ET_{i<j<n}x_i\ne x_j\land\forall y\,\Bigl(y\in z\eq\LOR_{i<n} y=x_i\Bigr)\Bigr)
\end{equation}
for all $n\in\omega$ is not weakly interpretable in~$\ec_L$, and consequently not interpretable in any
consistent existential theory.
\end{Cor}
\begin{Pf}
Apply Theorem~\ref{thm:main} to the formula interpreting $x\in z\land y\in z\land x\ne y$, and $R$ a chain longer than~$n$.
\end{Pf}

We can restate this in proper model-theoretic terminology. A theory~$T$ is said to \emph{eliminate $\exists^\infty$} (or
\emph{eliminate infinity}) if for every formula $\fii(\vc z,x)$, there exists $n$ such that for every
model $M\model T$ and $\vc a\in M$, if $\lh{\fii(\vc a,M)}\ge n$, then it is infinite. Specializing this to the
theory $T^\teq$, this means that for every $\fii(\vc z,\vc x)$ and $\psi(\vc x,\vc
y)$ (where $\Lh(\vc x)=\Lh(\vc y)=k$), there exists $n$ such that for every
$M\model T$ and $\vc a\in M$, if $\psi$ defines an equivalence relation on~$M^k$, and $\fii(\vc a,M^k)$ hits
at least $n$ equivalence classes, then it hits infinitely many.

As in Corollary~\ref{thm:existsinfty}, an application of Theorem~\ref{thm:main} to the formula $\fii(\vc z,\vc x)\land\fii(\vc
z,\vc y)\land\neg\psi(\vc x,\vc y)$ yields:
\begin{Cor}\label{cor:exinf}
$(\ec_L)^\teq$ has elimination of the $\exists^\infty$ quantifier.
\noproof\end{Cor}

Our second principal conclusion is the following tameness result on~$\ec_L$:
\begin{Cor}\label{thm:nsop3}
For any language~$L$, $\ec_L$ has $\prpty{NSOP_3}$. That is, the theory axiomatized by
\begin{gather*}
\forall x,y,z\,\neg(x<y\land y<z\land z<x),\\
\exists x_0,\dots,x_{n-1}\,\ET_{i<j<n}x_i<x_j
\end{gather*}
for $n\in\omega$ is not weakly interpretable in~$\ec_L$, and is not interpretable in any consistent existential theory.

Consequently, $\ec_L$ has the $\prpty{NSOP}$ property, i.e., no theory consistent with $\ec_L$ interprets a partial
order with arbitrarily long chains.
\end{Cor}
\begin{Pf}
Apply Theorem~\ref{thm:main} with $R$ being a directed $3$-cycle.
\end{Pf}

Let us state for the record that Corollaries \ref{thm:existsinfty} or~\ref{thm:nsop3} solve our original problem:
\begin{Cor}\label{cor:prf-r}
The theory $\tprf$ represents all partially recursive functions, but it does not interpret~$R$.
\noproof\end{Cor}

On the other hand, it should be stressed that $\ec_L$ is not \emph{that} tame, if the language $L$ is sufficiently
complicated (note that the observation below also stands in contrast to properties of random relational
structures, i.e., $\ec_L$ with $L$ purely relational, which is a simple theory, thus $\prpty{NTP_2}$).
\begin{Prop}\label{prop:tp2}
If $L$ contains an at least binary function symbol, then $\ec_L$ has~$\prpty{TP_2}$, hence it is not simple. More
generally, any theory weakly interpreting~$\trep\strf$ (i.e., with a consistent extension that loosely represents trf) has $\prpty{TP_2}$.
\end{Prop}
\begin{Pf}
It suffices to show the latter claim. Let $\num F(x,y)$ be the $\trep\strf$-function representing the recursive
function that interprets $x$ as a G\"odel number of a finite sequence, and outputs its $y$th element. Let $\vc a_{i,j}=(\num{i},\num{j})$, and $\fii(x,y_1,y_2)$ be the formula
\[\num F(x,y_1)=y_2.\]
Clearly, $\fii(x,\vc a_{i,j})\land\fii(x,\vc a_{i,k})$ is inconsistent for~$j\ne k$. On the other hand, if
$\sigma\in\omega^\omega$, and $n\in\omega$, let $s$ be the G\"odel number of~$\sigma\res n$. Then $\fii(\num s,\vc
a_{i,\sigma(i)})$ for all $i<n$. Thus, the type
\[\{\fii(x,\vc a_{i,\sigma(i)}):i<\omega\}\]
is consistent.
\end{Pf}

The assumption on~$L$ in Proposition~\ref{prop:tp2} is essential; see Theorem~\ref{thm:class} for more detailed model-theoretic
classification of the theories $\ec_L$ as $L$ varies.

Now, in the most general case when $L$ contains an at least
binary function symbol, there is still a gap left between Proposition~\ref{prop:tp2} and Corollary~\ref{thm:nsop3}. We can close it by
improving Corollaey~\ref{thm:nsop3} from $\prpty{NSOP_3}$ to $\prpty{NSOP_1}$, but the proof will no longer be self-contained: we will
rely on a characterization of $\prpty{NSOP_1}$ theories due to Chernikov and Ramsey~\cite{cher-ram} using an
independence relation in the spirit of the Kim--Pillay theorem.

We will work inside a \emph{monster model~$\monster$} of a completion $T\Sset\ec_L$, as in
Section~\ref{sec:model-theory}: a $\kappa$-saturated, strongly $\kappa$-homogeneous model of~$T$, where $\kappa$ is a
cardinal larger than $\|L\|$ and all structures we intend to handle; recall that subsets of~$\monster$ of size
$<\kappa$ are called small.
\begin{Def}\label{def:indep}
If $A$, $B$, and $C$ are small tuples (sequences), we say that \emph{$A$ is independent from $B$ over~$C$}, written as $A\indep[C]B$, if $\p{AC}\cap\p{BC}=\p C$, where $\p X$
denotes the substructure generated by~$X$, and the juxtaposition of two sequences denotes their concatenation. We will
often treat these tuples as sets where the context permits, seeing as the definition of $\indep$ does not depend on
their ordering.
\end{Def}

The definition of $A\indep[C]B$ is stated here in more general circumstances than what is required for~\cite{cher-ram}
(in particular, their characterization only needs the case when $C$ is a small \emph{model}, i.e., an elementary
submodel of the monster). We do it partly because we can---at no additional cost---and partly because we also want the
definition to conform to the shape of independence relations from the original Kim--Pillay theorem, which we will use
elsewhere in the paper. For the same reason, the next Lemma includes some properties of independence relations that are
not directly relevant to the characterization from~\cite{cher-ram}.
\begin{Lem}\label{lem:indep}
Let $\monster$ be a monster model of a completion $T$ of~$\ec_L$. The independence relation~$\indep$ has the following
properties for all small tuples $A,A',B,B',C,D$:
\begin{enumerate}
\item\label{item:41}
\paren{Invariance} If $f$ is an automorphism of~$\monster$, then $A\indep[C]B$ implies $f(A)\indep[{f(C)}]f(B)$.
\item\label{item:42}
\paren{Symmetry} $A\indep[C]B$ implies $B\indep[C]A$.
\item\label{item:43}
\paren{Monotonicity} If $A'\sset A$, and $B'\sset B$, then $A\indep[C]B$ implies $A'\indep[C]B'$.
\item\label{item:44}
\paren{Weak transitivity} $A\indep[B]C$ and $A\indep[BC]D$ implies $A\indep[B]CD$.
\item\label{item:45}
\paren{Existence} $A\indep[B]B$.
\item\label{item:46}
\paren{Strong finite character} If $A\nindep[C]B$, there is a formula $\fii(\vc x,\vc b,\vc c)\in\tp(A/BC)$ such
that $\vc a\nindep[C]B$ whenever $\monster\model\fii(\vc a,\vc b,\vc c)$.
\item\label{item:47}
\LP Extension%
\footnote{Our formulation of the extension property follows the statement of the Kim--Pillay
theorem~\cite{kim-pil,wagner}. As pointed out by the reviewer, this property is often postulated in a stronger form: for
any $A,B,C,D$, if $A\indep[C]B$, there is $A'\equiv_{BC}A$ such that $A'\indep[C]BD$. This easily follows
from~\ref{item:47} using \ref{item:41} and~\ref{item:44}.}%
\RP\ For any $A,B,C$, there is $A'\equiv_CA$ such that $A'\indep[C]B$.
\item\label{item:48}
\paren{Local character} For any $B$, and finite $A$, there is $B'\sset B$ such that $\lh{B'}\le\|L\|$, and
$A\indep[B']B$.
\item\label{item:49}
\paren{Independence theorem} If $A\indep[C]B$, $B\indep[C]B'$, $B'\indep[C]A'$, and $A'\equiv_CA$, there exists $A''$
such that $A''\equiv_{CB}A$, $A''\equiv_{CB'}A'$, and $A''\indep[C]BB'$.
\end{enumerate}
\end{Lem}
\begin{Pf}
Properties \ref{item:41}--\ref{item:45} are clear.

\ref{item:46}: By definition, $A\nindep[C]B$ implies that $t(\vc a,\vc c)=s(\vc b,\vc c)\notin\p C$ for some terms
$t,s$, and $\vc a\sset A$, $\vc b\sset B$, $\vc c\sset C$. Then we can take $t(\vc x,\vc c)=s(\vc b,\vc c)$ for
the formula $\fii$.

\ref{item:47}: We can extend the structure $\p{BC}$ with a disjoint copy of $\p{AC}\bez\p C$: that is, let us define a
structure~$D$ with domain
\[\p{BC}\cupd\{\ob x:x\in\p{AC}\bez\p C\},\]
with relations and functions defined so that they agree with the original structure on~$\p{BC}$, and so that $f=\id_{\p
C}\cup\ob x$ is an isomorphism of $\p{AC}$ to $\p C\cup\ob{\p{AC}\bez\p C}\sset D$. We can ensure $D$ is a
substructure of~$\monster$ (extending $\p{BC}$) using $\ec_L$ and $\kappa$-saturation. Then $A'=f(A)$ has the required
properties.

\ref{item:48}: Let us construct a chain $B_0\sset B_1\sset B_2\sset\cdots$ of subsets of~$B$ of size
$\lh{B_n}\le\lambda:=\|L\|$ as follows. We put $B_0=\nul$. Given $B_n$, let $M_n=\p{AB_n}\cap\p B$. Since
$\lh{M_n}\le\lambda$, there exists a set $B_n\sset B_{n+1}\sset B$ of size $\lh{B_{n+1}}\le\lambda$ such that
$M_n\sset\p{B_{n+1}}$.

Let $B'=\bigcup_{n\in\omega}B_n$. Then $\p{AB'}\cap\p B=\bigcup_nM_n\sset\p{B'}$ by construction, hence $A\indep[B']B$,
and $\lh{B'}\le\lambda$.

\ref{item:49}: In order to simplify the notation, we may assume without loss of generality that $C$, $A$, $A'$, $B$,
and $B'$ are structures, with $C\sset A,A',B,B'$. By the assumption, we have $A\cap B=A'\cap B'=B\cap B'=C$, and we can fix an isomorphism $f\colon A\simeq A'$ identical on~$C$. Put $D=\p{AA'BB'}$. We will extend $D$ into a model $D'$ with domain
\[D'=D\cupd\{\ob x:x\in A\bez C\}\cupd\{\ub x:x\in\p{AB}\bez(A\cup B)\}\cupd\{\ub{\ub x}:x\in\p{A'B'}\bez(A'\cup B')\}\]
using copies of parts of~$D$. We will also write $\ob{\ob y}$ for elements $y\in A'\bez C$, so that $\ob x=\ob{\ob y}$
if $f(x)=y$.

We define relations and functions on $D'$ so that
\begin{itemize}
\item $g=\id_B\cup\ob x\cup\ub x$ is an isomorphism of $\p{AB}$ to $B\cup\ob{A\bez C}\cup\ub{\p{AB}\bez(A\cup B)}$, and
\item $g'=\id_{B'}\cup\ob{\ob x}\cup\ub{\ub x}$ is an isomorphism of $\p{A'B'}$ to $B'\cup\ob{\ob{A'\bez C}}\cup\ub{\ub{\p{A'B'}\bez(A'\cup B')}}$.
\end{itemize}
It is important to note there is no conflict between the two clauses: the intersection of the two targets is
$C\cup\ob{A\bez C}=C\cup\ob{\ob{A'\bez C}}$, which is asked to be made isomorphic to $A$ via $\id_C\cup\ob x$, and to
$A'$ via $\id_C\cup\ob{\ob x}=\id_C\cup\ob{f^{-1}(x)}$; these two requirements are equivalent, as the two mappings
commute with the isomorphism $f\colon A\to A'$.

Now, using $\ec_L$ and $\kappa$-saturation, we can embed $D'$ as a substructure of~$\monster$ extending~$D$. Let
$A''=g(A)=g'(A')$. Then $g$ is an isomorphism of $A$ to~$A''$ identical on $B\Sset C$, thus $A''\equiv_{CB}A$, and
similarly $A''\equiv_{CB'}A'$ via~$g'$. Finally, $\p{A''C}=A''$, $\p{BB'C}\sset D$, and $A''\cap D=C$, thus
$A''\indep[C]BB'$.
\end{Pf}

The following is a restatement of Proposition~5.8 in Chernikov and Ramsey~\cite{cher-ram}.
\begin{Thm}\label{thm:cher-ram}
Let $\monster$ be a monster model of a complete theory~$T$, and $A\indep[M]B$ an independence relation on small tuples
$A,B$, and small models $M\model T$, that satisfies the appropriate restrictions of properties \ref{item:41},
\ref{item:42}, \ref{item:43}, \ref{item:45}, \ref{item:46}, and \ref{item:49} from Lemma~\ref{lem:indep}. Then $T$ is
$\prpty{NSOP_1}$.
\noproof\end{Thm}
\begin{Cor}\label{cor:nsop1}
For any language~$L$, $\ec_L$ is $\prpty{NSOP_1}$.
\noproof\end{Cor}
\begin{Rem}\label{rem:kim-indep}
Using Theorem~9.1 in~\cite{kap-ram}, it can be seen that our $\indep$ coincides with the relation of \emph{Kim-independence}.
\end{Rem}

We mention that Corollary~\ref{cor:nsop1} was independently discovered by Kruckman and Ramsey~\cite{kruck-ram}, who learned of
the problem from an earlier unpublished version of this paper where it was posed as an open problem.

\section{Conclusion}\label{sec:conclusion}

We succeeded in our original goal of separating interpretability of~$R$ from representability of recursive functions.
More generally, we obtained a criterion for interpretability of $\exists\forall$ theories in existential theories,
showing in particular that we may assume such interpretations to be quantifier-free. We believe these results are
interesting in their own right, of course, but at the same time we place as much value on the connection between formal
arithmetic and model theory that it revealed: while model-theoretic methods are often used in the study of arithmetic,
typically this means to work with models of the (fairly strong) theories of arithmetic themselves, which are quite
unlike the kind of tame model theory we encountered in this paper. It would be interesting to see if more such
connections are waiting to be discovered.

\appendix
\section{Elimination of imaginaries}\label{sec:elim-imag}

In this section, we discuss elimination of imaginaries in the theories~$\ec_L$. We put it here in the appendix as it is
rather tangential to our main topic; we only need it in the proof of Theorem~\ref{thm:ea}.

One way to describe elimination
of imaginaries is that, loosely speaking, it allows to replace any interpretation with an interpretation with absolute
equality. We recall the proper definition below, along with some important variants of the notion.
\begin{Def}\label{def:ei}
A theory $T$ has \emph{elimination of imaginaries (e.i.)}  if for every $M\model T$ and $e\in M^\teq$,
there is a tuple $\vc b\in M$ such that $\dcl_{M^\teq}(e)=\dcl_{M^\teq}(\vc b)$.

$T$ has \emph{weak e.i.}\  if for every $M\model T$ and $e\in M^\teq$,
there is $\vc b\in M$ such that $\vc b\in\acl_{M^\teq}(e)$ and $e\in\dcl_{M^\teq}(\vc b)$.

As a special case, for any $k,l>0$ let $\sim_{k,l}$ be the equivalence relation on injective $l$-tuples of
$k$-tuples (represented as $kl$-tuples) defined by
\[\p{\vc a_i:i<l}\sim_{k,l}\p{\vc b_i:i<l}\iff\{\vc a_i:i<l\}=\{\vc b_i:i<l\},\]
so that $M^{kl}/{\sim_{k,l}}$ represents $l$-element subsets of~$M^k$. If for all $M\model T$ and all $k,l>0$, every
$a\in M^{kl}/{\sim_{k,l}}$ is interdefinable (in $M^\teq$) with a tuple $\vc b\in M$, then $T$ has \emph{coding of
finite sets}. See~\cite{cas-far} for an exposition of various forms of e.i.; in particular, $T$ has e.i.\ iff it has
weak e.i.\ and coding of finite sets.
\end{Def}

It would be nice if $\ec_L$ had e.i., however this is too good to be true: as we will prove shortly, the theory does not
have coding of finite sets.
\begin{Lem}\label{lem:acl}
If $M\model\ec_L$, and $A\sset M$, then $\acl_M(A)=\p A$.
\end{Lem}
\begin{Pf}
Assume $A$ is a submodel, and $M\model\fii(b,\vc a)$, where $b\in M\bez A$, and $\fii$ is open. Let
$N=M\cupd(\omega\times(M\bez A))$, where each $A\cup(\{i\}\times(M\bez A))$ is an isomorphic copy of~$M$, as in the
proof of Lemma~\ref{lem:indep}~\ref{item:47}. Then $N\model\fii(\p{i,b},\vc a)$ for each $i<\omega$, hence
$\fii(M,\vc a)$ is also infinite by existential closedness.
\end{Pf}
\begin{Lem}\label{lem:indi}
Every $M\model\ec_L$ has elementary extensions with arbitrarily large sets~$X$ totally indiscernible over~$M$.
Moreover, we can choose $X$ to have the additional property that any $L_M$-term $t(x_1,\dots,x_k)$ containing all the
indicated variables defines an injective function on~$X^k$.
\end{Lem}
\begin{Pf}
Let $M\model\ec_L$, and $X$ be a set disjoint from~$M$. An $L_M$-term in variables $X$ is \emph{reduced} if
it has no constant subterms other than $M$-constants; that is, variables and $M$-constants are reduced, and
if $t_1,\dots,t_k$ are reduced terms, and $F$ is a $k$-ary function symbol of~$L$, then $F(t_1,\dots,t_k)$ is reduced unless
all the $t_i$ are $M$-constants.

Let $R$ be the model whose domain is the set of all reduced terms, with realizations of relations the same as in~$M$ (i.e.,
unsatisfied by tuples involving any non-constant terms), and functions realized in the obvious way. Let $N$ be an e.c.\
extension of~$R$. Since every permutation of~$X$ extends to an automorphism of~$R$ fixing~$M$, any tuples of distinct elements
$x_1,\dots,x_k$ and $x'_1,\dots,x'_k$ of~$X$ satisfy the same atomic formulas with parameters from~$M$. Thus in~$N$,
$X$ is a totally indiscernible set over~$M$ by quantifier elimination.
\end{Pf}
\begin{Prop}\label{prop:no-ei}
Neither $\ec_L$ nor any its completion has coding of unordered pairs of elements, and a fortiori elimination of
imaginaries.
\end{Prop}
\begin{Pf}
Let $N$ be a model of~$\ec_L$ with $\{a,b\}$ a $2$-element totally indiscernible set satisfying the property from Lemma~\ref{lem:indi}. Without loss of generality, $N$ is strongly $\omega$-homogeneous.
Assume for contradiction that there is a tuple $\vc u$ in~$N$
interdefinable with the representation of $\{a,b\}$ (that is, with the element $\p{a,b}/{\sim_{1,2}}$ of~$N^\teq$). By
Lemma~\ref{lem:acl}, all elements of~$\vc u$ are in the submodel generated by~$a,b$, hence they are given by a tuple of
terms $\vc t(a,b)$. By $a,b\equiv b,a$ and homogeneity, there is an automorphism $f$ such that $f(a)=b$ and $f(b)=a$;
since $f$ preserves $\{a,b\}$, it also preserves $\vc u$, hence $t_i(a,b)=t_i(b,a)$. Using the extra property, this
can only happen if all the $t_i$ are closed terms. Thus $\vc u$, hence $\{a,b\}$, is $\nul$-definable, and
$a,b\in\acl_N(\nul)$; using Lemma~\ref{lem:acl} again, $a$ and~$b$ are in fact values of closed terms, but this
contradicts $a\equiv b$.
\end{Pf}

Short of full e.i., the next best thing we can hope for is weak e.i. This will turn out to hold for~$\ec_L$, and
thankfully it is still enough for our intended application.

Let us work again in a monster model~$\monster$ of a completion of~$\ec_L$. Recall the following
characterization~\cite[Facts~1.2]{cas-far}: weak e.i.\ holds iff for every relation~$R$ definable with parameters
in~$\monster$, there exists a smallest algebraically closed set defining~$R$. (By Lemma~\ref{lem:acl}, algebraically
closed set = substructure for us.) We observe easily that $\ec_L$ satisfies a somewhat weaker property:
\begin{Lem}\label{lem:def-directed}
For any definable relation~$R(\vc x)$, the class of substructures that define~$R$ is directed. That is, if $R$ is
definable over~$\vc b$, and over~$\vc c$, it is also definable over $\p{\vc b}\cap\p{\vc c}$.
\end{Lem}
\begin{Pf}
In order to simplify the notation, we will omit bars over finite tuples.

Put $A=\p b\cap\p c$. By saturation, it suffices to show that
\begin{equation}\label{eq:4}
x\equiv_Ax'\implies\bigl(R(x)\eq R(x')\bigr).
\end{equation}
Assume that $R(x)\eq\fii(b,x)\eq\psi(c,x)$.
\pagebreak[2]
\begin{Cl}
\ \begin{enumerate}
\item\label{item:50}
If $b\equiv_cb'$, then $R(x)\eq\fii(b',x)$.
\item\label{item:53}
If $c\equiv_bc'$, and $b\equiv_{c'}b'$, then $R(x)\eq\fii(b',x)$.
\item\label{item:51}
If $b\equiv_Ab'$, and $x\indep[A]b'$, then $R(x)\eq\fii(b',x)$.
\item\label{item:52}
If $x\equiv_Ax'$, and $x\indep[A]b$, then $R(x)\eq R(x')$.
\end{enumerate}
\end{Cl}
\begin{Pf*}
\ref{item:50}: Let $\alpha(v,w)$ denote $\forall x\,(\fii(v,x)\eq\psi(w,x))$. Since $\alpha(b,c)$, and $b\equiv_cb'$, we have
$\alpha(b',c)$.

\ref{item:53}: By the dual statement to~\ref{item:50}, we have $R(x)\eq\psi(c',x)$. Thus, $\alpha(b,c')$, whence
$\alpha(b',c')$.

\ref{item:51}: By Lemma~\ref{lem:indep}~\ref{item:47}, there exists $c'\equiv_bc$ such that $c'\indep[b]x$. Since
$c\indep[A]b$ by the definition of~$A$, $c'\equiv_bc$ implies $c'\indep[A]b$. By Lemma~\ref{lem:indep}~\ref{item:44}, we
obtain $c'\indep[A]x$. Thus, by Lemma~\ref{lem:indep}~\ref{item:49}, there exists $b''$ such that $b''\equiv_{c'}b$, and
$b''\equiv_{Ax}b'$. Using \ref{item:53}, this implies
\[R(x)\eq\fii(b'',x)\eq\fii(b',x).\]

\ref{item:52}: There exists $b'$ such that $xb\equiv_Ax'b'$. Then $x'\indep[A]b'$, thus
\[R(x)\eq\fii(b,x)\eq\fii(b',x')\eq R(x')\]
by~\ref{item:51}.
\end{Pf*}

Finally, to prove~\eqref{eq:4}, let $x''\equiv_Ax$ be such that $x''\indep_Ab$ using
Lemma~\ref{lem:indep}~\ref{item:47}. Then $R(x)\eq R(x'')\eq R(x')$ by \ref{item:52} of the Claim.
\end{Pf}

Incidentally, the previous lemma implies another property: 
$T$ has \emph{Galois e.i.}\ if for every $M\model T$, $\vc a\in M$, and $e\in\acl_{M^\teq}(\vc a)$,
there is $\vc b\in M$ such that $\dcl_{M^\teq}(\vc ae)=\dcl_{M^\teq}(\vc a\vc b)$.

\begin{Prop}\label{prop:gal-ei}
For any language~$L$, $\ec_L$ has Galois elimination of imaginaries.
\end{Prop}
\begin{Pf}
By \cite[Prop.~3.9]{cas-far}, Galois e.i.\ is equivalent to the conjunction of elimination of strong types (ST), and
coding of Galois finite sets. The latter follows from our Lemma~\ref{lem:acl} by~\cite[Facts~3.7]{cas-far}. It thus
suffices to show ST; by \cite[Prop.~3.2]{cas-far}, this is equivalent (in view of Lemma~\ref{lem:acl}) to the claim
that $\acl_{M^\teq}(A)=\dcl_{M^\teq}(A)$ for real sets~$A$ in a monster model~$\monster\model\ec_L$. We can reformulate this
as follows: if $E$ is an equivalence relation (on $k$-tuples) definable over~$A$ with finitely many equivalence
classes, then all the equivalence class are individually definable over~$A$.

So, let us fix a $k$-tuple $\ob x$, we will show that the equivalence class of~$\ob x$ is $A$-definable. Assume that
$E$ has $n$~classes. By repeated use of Lemma~\ref{lem:indep}~\ref{item:47}, we can find a sequence $\{\ob x^i:i\le n\}$
of tuples such that $\ob x\equiv_A\ob x^i$, and $\ob x^i\indep[A]\ob x^0\dots\ob x^{i-1}\ob x$. By the pigeonhole
principle, there exist $0\le i<j\le n$ such that $\ob x^i$ and $\ob x^j$ are in the same class~$C$ of~$E$. It follows
that $C$ is definable over $A\ob x^i$ and over $A\ob x^j$, hence by Lemma~\ref{lem:def-directed}, it is definable over
$\p{A\ob x^i}\cap\p{A\ob x^j}=\p A$, i.e., over~$A$. But then $\ob x\equiv_A\ob x^i$ implies that $\ob x\in C$.
\end{Pf}

It may come as anticlimactic that we will state the result we are most interested in, viz.\ weak e.i.\ for~$\ec_L$, without proof: while the author has
figured out a long and cumbersome argument, it was independently shown in an easier way by Kruckman and
Ramsey~\cite[\S3.4]{kruck-ram}, and we invite the interested reader to consult their paper.
\begin{Thm}\label{thm:wei}
For any language~$L$, $\ec_L$ has weak elimination of imaginaries.
\noproof\end{Thm}

We end this section by stating explicitly the consequence of Theorem~\ref{thm:wei} for definable functions of~$\ec_L^\teq$
that we will need in the proof of Theorem~\ref{thm:ea}. While it is somewhat hairy to formulate, it follows by a simple
compactness argument.
\begin{Lem}\label{lem:terms}
Let $T\Sset\ec_L$. If 
\[T\vdash\alpha(\vc x)\to\exists^{=m}\vc u\,\beta(\vc u,\vc x),\]
there are formulas $\{\alpha_i(\vc x):i<n\}$
for some $n>0$, and for each $i<n$, tuples of terms $\vc t^{i,j}(\vc x)$, $j<m$,
such that $T$ proves that $\{\alpha_i:i<n\}$ define a partition of $\alpha$ \paren{i.e., $\alpha(\vc
x)\eq\LOR_{i<n}\alpha_i(\vc x)$, and $\alpha_i(\vc x)\to\neg\alpha_{i'}(\vc x)$ for $i\ne i'$},
and
\[\alpha_i(\vc x)\to\Bigl(\beta(\vc u,\vc x)\eq\LOR_{j<m}\vc u=\vc t^{i,j}(\vc x)\Bigr).\]
Moreover, $T$ proves $\alpha_i(\vc x)\to\vc t^{i,j}(\vc x)\ne\vc t^{i,j'}(\vc x)$ for $j\ne j'$.
\end{Lem}
\begin{Pf}
Note that $\alpha_i(\vc x)\to\vc t^{i,j}(\vc x)\ne\vc t^{i,j'}(\vc x)$ for $j\ne j'$ follows from the rest,
specifically $\alpha(\vc x)$ implies that there are \emph{at least} $m$~tuples satisfying $\beta(\vc u,\vc x)$.

We will prove the statement by induction on~$m$. If $m=0$, there is nothing to prove. Assuming the result holds
for~$m$, we will prove it for $m+1$.

By Lemma~\ref{lem:acl} and a
compactness argument, there are tuples of terms $\vc t^{i,m}(\vc x)$, $i<n$, such that $T$ proves
\[\alpha(\vc x)\to\LOR_{i<n}\beta(\vc t^{i,m}(\vc x),\vc x).\]
Putting
\[\alpha_i(\vc x)\eq\beta(\vc t^{i,m}(\vc x),\vc x)\land\ET_{i'<i}\neg\beta(\vc t^{i',m}(\vc x),\vc x),\]
we have a partition of $\alpha$ into formulas $\alpha_i$ such that $T$ proves
\[\alpha_i(\vc x)\to\beta(\vc t^i(\vc x),\vc x).\]
Let us write
\[\beta_i(\vc u,\vc x)\eq\beta(\vc u,\vc x)\land\vc u\ne\vc t^i(\vc x).\]
Then $T$ proves
\[\alpha_i(\vc x)\to\exists^{=m}\vc u\,\beta_i(\vc u,\vc x).\]
Using the induction hypothesis, we can further refine the partition so that there are tuples of terms $\vc
t^{i,j}(\vc x)$, $j<m$, such that
\[T\vdash\alpha_i(\vc x)\to\Bigl(\beta_i(\vc u,\vc x)\eq\LOR_{j<m}\vc u=\vc t^{i,j}(\vc x)\Bigr).\]
Then $T$ proves
\[\alpha_i(\vc x)\to\Bigl(\beta(\vc u,\vc x)\eq\LOR_{j<m+1}\vc u=\vc t^{i,j}(\vc x)\Bigr)\]
as required.
\end{Pf}

\begin{Prop}\label{prop:equi-expl}
Let $T\Sset\ec_L$, and $E$ be an equivalence relation on $k$-tuples definable in~$T$. Then there are
\begin{itemize}
\item integers $n>0$, $r>0$, and $m_i\ge0$ for $i<n$,
\item formulas $\fii_i(\vc x)$ for $i<n$, and
\item terms $t^{i,j}_l(\vc x)$ for $i<n$, $j<m_i$, $l<r$,
\end{itemize}
such that $T$ proves
\begin{equation}\label{eq:5}
E(\vc x,\vc y)\eq\LOR_{i<n}\Bigl(\fii_i(\vc x)\land\fii_i(\vc y)\land
    \bigl\{\vc t^{i,j}(\vc x):j<m_i\bigr\}=\bigl\{\vc t^{i,j}(\vc y):j<m_i\bigr\}\Bigr),
\end{equation}
where $\vc t^{i,j}(\vc x)$ denotes the tuple $\p{t^{i,j}_l(\vc x):l<r}$, and the following condition hold:
\begin{itemize}
\item The formulas $\fii_i$ form a partition, i.e., $T$ proves $\LOR_{i<n}\fii_i(\vc x)$, and $\fii_i(\vc
x)\to\neg\fii_{i'}(\vc x)$ for $i\ne i'$.
\item On each part $\fii_i$, the tuples $\vc t^{i,j}$ are pairwise distinct, i.e.,
\[T\vdash\fii_i(\vc x)\to\LOR_{l<r}t^{i,j}_l(\vc x)\ne t^{i,j'}_l(\vc x)\]
for each $i<n$, and $j<j'<m_i$.
\end{itemize}
\end{Prop}
\begin{Pf}
By Theorem~\ref{thm:wei}, for each $M\model T$ and an equivalence class $e$ of~$E$, there is a real tuple $\vc a\in M$ such that
$\vc a\in\acl_{M^\teq}(e)$, and $e\in\dcl_{M^\teq}(\vc a)$. Thus, for each such $M$, $e$, and $\vc a$, we can find a
formula $\psi(\vc u,\vc x)$ and an integer $m\ge1$ such that:
\begin{itemize}
\item $\psi(\vc a,\vc x)$ for all $\vc x\in e$.
\item $\psi(\vc u,\vc x)\to\bigl(\psi(\vc u,\vc y)\eq E(\vc x,\vc y)\bigr)$; that is, for a given $\vc u$, the
set $\{\vc x:M\model\psi(\vc u,\vc x)\}$ is either empty, or an equivalence class of $E$.
\item For a given $\vc x$, there are either none or exactly $m$ tuples $\vc u$ such that $\psi(\vc u,\vc x)$.
\end{itemize}
Using a compactness argument, there is a finite
partition coarser than~$E$ definable in $T$ by formulas $\{\fii_i(\vc x):i<n\}$, and for each $i<n$, there is a formula
$\psi_i(\vc u,\vc x)$, and an integer $m_i\ge1$, such that $T$ proves:
\begin{gather*}
\psi_i(\vc u,\vc x)\to\bigl(\psi_i(\vc u,\vc y)\eq E(\vc x,\vc y)\bigr),\\
\psi_i(\vc u,\vc x)\to\fii_i(\vc x),\\
\fii_i(\vc x)\to\exists^{=m_i}\vc u\,\psi_i(\vc u,\vc x).
\end{gather*}
For a given~$i$, if the tuple $\vc u$ in $\psi_i(\vc x,\vc u)$ has length~$0$, then $\fii_i(\vc x)$ is equivalent
to~$\psi_i(\vc x)$, and it defines a single equivalence class of~$E$; in this case, we formally replace $\psi_i$ with an always false formula using a dummy variable~$u$,
and put $m_i=0$. In this way, we may ensure that all the $\vc u$ are nonempty; by repeating one of its elements, we
may in fact assume that they all have the same length $l>0$ independent of~$i$.

Using Lemma~\ref{lem:terms}, after possibly refining the partition $\{\fii_i:i<n\}$, we can find tuples of terms $\vc
t^{i,j}(\vc x)$ such that $T$ proves
\[\psi_i(\vc u,\vc x)\eq\fii_i(\vc x)\land\LOR_{j<m_i}\vc u=\vc t^{i,j}(\vc x).\]
This implies~\eqref{eq:5}, that is,
\[\fii_i(\vc x)\land\fii_i(\vc y)\to\Bigl(E(\vc x,\vc y)\eq\bigl\{\vc t^{i,j}(\vc x):j<m_i\bigr\}=\bigl\{\vc
t^{i,j}(\vc y):j<m_i\bigr\}\Bigr).\]
Indeed, $\{\vc t^{i,j}(\vc x):j<m_i\}=\{\vc u:\psi_i(\vc u,\vc x)\}$, and similarly for $\vc y$. The properties
of $\psi_i$ ensure that if $E(\vc x,\vc y)$, these sets are equal, whereas if not (which can only happen if $m_i>0$),
they are disjoint, hence distinct.
\end{Pf}

\section{Dependence on language}\label{sec:dependence-language}

The results in Section~\ref{sec:classification-ec_l} give fairly tight model-theoretic classification of $\ec_L$
($\prpty{NSOP_1}$, but $\prpty{TP_2}$) in the
case that $L$ includes at least one at least binary function symbol. However, the theories behave in different ways for
other languages~$L$; in this section, we summarize the main model-theoretic properties of the theories $\ec_L$ in
dependence on~$L$.

For the benefit of readers coming from a non-model-theoretic background, we recall that a theory~$T$ is
\emph{$\kappa$-stable} if for every $M\model T$ and $A\sset M$ of size $\lh A\le\kappa$, there are at most $\kappa$ complete types over~$A$. We say that $T$ is \emph{stable} if it is $\kappa$-stable for some infinite
cardinal~$\kappa$, and \emph{superstable} if it is $\kappa$-stable for all sufficiently large cardinals~$\kappa$. An
even stronger condition is that $T$ be \emph{totally transcendental}; officially, this means that every formula~$\fii$
has Morley rank $\mathrm{MR}(\fii)<\infty$, but we can use the
following characterization: $T$ is totally transcendental iff all countable-language fragments of~$T$ are
$\omega$-stable (which implies $\kappa$-stable for all $\kappa\ge\omega$).
The class of $\omega$-stable theories includes uncountably
categorical theories, which in turn include \emph{strongly minimal} theories (meaning that for all $M\model T$, the
only subsets of~$M$ definable with parameters are finite or cofinite).

\begin{Thm}\label{thm:class}
Let $L$ be a language. Then any complete extension of~$\ec_L$ is
\begin{enumerate}
\item\label{item:19}
strongly minimal
iff $L$ consists of nullary symbols;
\item\label{item:20}
totally transcendental iff $L$ consists of nullary symbols, and either finitely many unary relations, or one unary
function;
\item\label{item:21}
superstable iff $L$ consists of at most unary symbols, at most one of which is a unary function;
\item\label{item:22}
stable iff it has $\prpty{NIP}$ iff $L$ consists of at most unary symbols;
\item\label{item:24}
supersimple iff $L$ consists of relations, constants, and at most one unary function;
\item\label{item:23}
simple iff it has $\prpty{NTP_2}$ iff $L$ consists of relations and at most unary functions;
\item\label{item:25}
$\prpty{NSOP_1}$.
\end{enumerate}
\end{Thm}
\begin{Pf}
We may assume $L$ contains no nullary relations, as these are fixed to true or false in any complete extension.

First, we establish the right-to-left implications. \ref{item:19} is obvious. \ref{item:20}: We may assume $L$ is
countable, we will verify $\ec_L$ is $\omega$-stable. Let $M\model\ec_L$, and $A\sset M$ be a set of parameters, which
we may assume to be a submodel. If $L$ consists of constants and finitely many unary relations~$P_i(x)$, a type
$\tp(b/A)$ of an element $b\notin A$ is determined by $\{i:P_i(b)\}$, hence there are only finitely many (plus $\lh A$
trivial types for $b\in A$). If $L$ consists of constants and a unary function $F(x)$, $\tp(b/A)$ is determined by the
least $n$ such that $F^n(b)\in A$ and the value of $F^n(b)$, or in case there is no such~$n$, by the shape of the chain
$\{F^n(b):n<\omega\}$. This makes $\lh A+\aleph_0$ possibilities.

\ref{item:21}\footnote{A theory is superstable iff it is stable and supersimple, hence \ref{item:21} follows from
\ref{item:22} and~\ref{item:24}. However, we prefer to give a direct proof not relying on sophisticated tools like the
Kim--Pillay theorem.}: Assume $L$ consists of constants, unary relations $\{P_i(x):i\in I\}$, and one unary function
$F(x)$. Let $A\sset M\model\ec_L$ be a submodel, and $b\in M$. Similarly to the totally transcendental case, $\tp(b/A)$ is
determined by the type of $b$ over~$A$ in the relation-free reduct of~$M$, for which there are $\lh A+\aleph_0$
possibilities, and by $\{\p{i,n}:i\in I,n<\omega,P_i(F^n(b))\}$. Thus, $\lh{S_1(A)}\le\lh A+2^{\|L\|}$.

\ref{item:22}: Let $L$ be unary, and $M,A,b$ as before. Atomic formulas involved in~$\tp(b/A)$ are of the forms
$t(b)=s(b)$ or $P(t(b))$ (i.e., not referring to~$A$, hence at most $\|L\|$ many), or $t(b)=a$ for $a\in A$. For any $t(x)$,
the type either contains one formula of the form $t(b)=a$, or all of $\{t(b)\ne a:a\in A\}$. Thus, $\lh{S_1(A)}\le\lh
A^{\|L\|}$; in particular, $\ec_L$ is $\kappa$-stable whenever $\kappa=\kappa^{\|L\|}$.

\ref{item:23}: The Kim--Pillay theorem~\cite{kim-pil} (see also \cite[2.6.1]{wagner}) states that a theory is simple if
we can define an independence relation $A\indep[C]B$ that satisfies properties \ref{item:41}, \ref{item:42},
\ref{item:44}, [a weaker form of] \ref{item:46}, \ref{item:47}, \ref{item:48}, and \ref{item:49} from
Lemma~\ref{lem:indep}, and the converse implication to \ref{item:44}; the latter amounts to \emph{base
monotonicity}: if $C'\sset B$, then $A\indep[C]B$ implies $A\indep[CC']B$.

Now, if $L$ contains no functions of arity $2$~or more, the relation $A\indep[C]B$ from Definition~\ref{def:indep} can be
restated as
\[A\indep[C]B\iff\forall a\in A\,\forall b\in B\,\p a\cap\p b\sset\p C.\]
This shows that it satisfies base monotonicity, even in the stronger form that $A\indep[C]B$ implies $A\indep[CC']B$ for
arbitrary~$C'$.

\ref{item:24}: In terms of the independence relation, a simple theory is supersimple iff it satisfies a strong form of
local character: for every $B$ and finite~$A$, there is
a \emph{finite} $B'\sset B$ such that $A\indep[B']B$ (cf.~\cite{wagner}). If $L$ consists of relations and constants,
we can take $B'=A\cap B$. If $L$
also contains one unary function~$F$, we construct a set $B'\sset B$ with $\lh{B'}\le\lh A$ such that for each
$u\in A$, if $F^n(u)\in\p B\bez\p\nul$ for some~$n$, then the least such~$n$ satisfies $F^n(u)=F^m(v)$ for
some~$v\in B'$ and~$m\in\omega$.

\ref{item:25} is Corollary~\ref{cor:nsop1}.

Now we turn to the left-to-right implications.

\ref{item:23} follows from Proposition~\ref{prop:tp2}.

\ref{item:24}: We may assume $T$ is simple, i.e., all functions in~$L$ are at most unary. If there are two unary
functions $F(x)$, $G(x)$, a monster model $\monster\model T$ will contain an $a$ such that $t(a)\ne s(a)$ for any pair of
distinct unary terms $t$, $s$. Let $B=\{G(F^n(a)):n<\omega\}$. Then any finite $C\sset B$ satisfies $\p C\cap B=C\ssset
B\sset\p a$, hence $a\nindep[C]B$. Thus $T$ is not supersimple.

\ref{item:22}: If $L$ contains an at least binary function, it is even $\prpty{TP_2}$. If $L$ contains an at least
binary relation, wlog $P(x,y)$, then the formula $P(x,y)$ has $\prpty{IP}$, as every finite model embeds into any model
of~$\ec_L$.

\ref{item:21}: We may assume $L$ is unary, lest the theory is not even stable. So, assume $L$ contains two unary
functions $F(x)$, $G(x)$. Let $\kappa$ be an arbitrarily large cardinal such that $\kappa^\omega>\kappa$,
$M\model\ec_L$, and $A=\{a_\alpha:\alpha<\kappa\}\sset M$ of cardinality~$\kappa$. For every
$\sigma\colon\omega\to\kappa$, the type
\[p_\sigma(x)=\{G(F^n(x))=a_{\sigma(n)}:n\in\omega\}\]
over~$A$ is consistent, and $p_\sigma$ and $p_\tau$ are incompatible for $\sigma\ne\tau$. Thus,
$\lh{S_1(A)}\ge\kappa^\omega$, and $T$ is not $\kappa$-stable.

\ref{item:20}: We may assume $T$ is superstable. If $L$ contains a function $F(x)$, and a relation $P(x)$, the
$2^\omega$ types
\[p_I(x)=\{P(F^n(x)):n\in I\}\cup\{\neg P(F^n(x)):n\notin I\}\]
for $I\sset\omega$ witness that the $\p{F,P}$-fragment of~$T$ is not $\omega$-stable, hence $T$ is not totally
transcendental. If $L$ contains infinitely many predicates $P_n(x)$, $n<\omega$, we can likewise use the types
$\{P_n(x):n\in I\}\cup\{\neg P_n(x):n\notin I\}$.

\ref{item:19}: If $L$ contains a non-nullary (wlog unary) predicate $P(x)$, then both $P(M)$ and its complement are
infinite for any $M\model\ec_L$, thus $M$ is not minimal. Likewise, if $L$ contains a function $F(x)$, the formula
$F(x)=x$ defines an infinite set with infinite complement.
\end{Pf}

Apart from tameness properties from classification theory, we also discuss some more elementary invariants of the
theories, namely the number of types, and the number of complete extensions. We only consider countable languages for
the rest of this section.

Recall that a countable complete theory is called \emph{small} if it has countable many complete $n$-types for all
$n\in\omega$; this holds if and only if it has a saturated countable model.
\begin{Prop}\label{prop:categ}
Let $L$ be countable, and $T$ a complete extension of~$\ec_L$.
\begin{enumerate}
\item\label{item:26}
Let $L$ consist of relations and constants, where the number of non-nullary relations and $T$-unequal constants is
finite. Then $T$ is $\omega$-categorical.
\item\label{item:27}
Let $L$ consist either of nullary symbols and one unary function, or of nullary relations, finitely many unary
relations, and infinitely many $T$-unequal constants. Then $T$ has $\aleph_0$ complete $n$-types for each
$0<n<\omega$, hence it is not $\omega$-categorical, but it is small.
\item\label{item:28}
Otherwise $T$ has $2^\omega$ complete $1$-types.
\end{enumerate}
\end{Prop}
\begin{Pf}

\ref{item:26}: By quantifier elimination, there are only finitely many formulas in $n$~variables for every~$n$
(ignoring sentences).

\ref{item:27}: On the one hand, if $F$ is a unary function, there are infinitely many incompatible $1$-types of the
form
\[F^n(x)=F^{n+1}(x)\land\ET_{i<j\le n}F^i(x)\ne F^j(x).\]
If $\{c_i:i<\omega\}$ are provably pairwise distinct constants, then we have infinitely many $1$-types extending
$x=c_i$.

On the other hand, $T$ is $\omega$-stable by Theorem~\ref{thm:class}.

\ref{item:28}: If we have a unary function $F(x)$, and a nonnullary relation (wlog unary) $R(x)$, we have countably
many independent atomic formulas $R(F^n(x))$, thus $2^\omega$ $1$-types. Likewise, if there is another unary
function $G(x)$, the formulas $G(F^n(x))=x$ are independent, and if we have a binary (or more) function
$H(x,y)$, we can consider the formulas $H(s_1(x),s_n(x))=x$, where $s_1(x):=H(x,x)$, and $s_{n+1}(x):=H(x,s_n(x))$.
If there are infinitely many proper predicates, or an at least binary predicate and infinitely many distinct constants,
we are also done.
\end{Pf}

\begin{Cor}\label{prop:complext}
Let $L$ be countable.
\begin{enumerate}
\item\label{item:29}
$\ec_L$ is complete iff $L$ contains no nullary symbols, or consists of one constant.
\item\label{item:30}
If $L$ contains no constants, and only finitely many nullary relation symbols (with no restrictions on non-nullary
symbols), then $\ec_L$ has finitely many complete extensions.
\item\label{item:15}
If $L$ is finite, and contains no nonconstant functions, then $\ec_L$ has finitely many complete
extensions.
\item\label{item:31}
If $L$ consists of a unary function, and finitely many nullary symbols, at least one of which is a constant, then
$\ec_L$ has countably infinitely many complete extensions.
\item\label{item:32}
Otherwise $\ec_L$ has $2^\omega$ complete extensions.
\end{enumerate}
\end{Cor}
\begin{Pf}
\ref{item:29}: It is easy to see that in the other cases, there is at least one nontrivial atomic sentence.

\ref{item:30}--\ref{item:32}: If there are infinitely many nullary relations or constants, there are infinitely many
independent atomic sentences, hence $2^\omega$ complete extensions.

Assume there are only finitely many nullary
symbols, and let $L_0$ be $L$ minus constants. The only atomic sentences in~$L_0$ are nullary
relations, hence $\ec_{L_0}$ has finitely many complete extensions. If $L$ has $k>0$ constants, then completions
of~$\ec_L$ correspond to complete $k$-types over completions of~$\ec_{L_0}$. By Proposition~\ref{prop:categ}, these are finitely
many if $L_0$ is finite and contains no functions, $\aleph_0$ if $L_0$ consists of one unary function and nullary
predicates, and $2^\omega$ otherwise.
\end{Pf}

Let us touch upon a somewhat different topic now. In Section~\ref{sec:model-compl-empty}, we proved the existence of the model completion~$\ec_L$ in a
laborious way by, essentially, computing an explicit finite axiomatization of resultants. Model theorists are not very
keen on getting their hands dirty with actual formulas, and prefer higher-level methods; in particular, a very popular
technique for construction of model completions is using \emph{Fra\"\i ss\'e limits}. For example, if $L$ is a finite
relational language, one can show the existence of $\ec_L$ quite easily by taking the theory of the Fra\"\i ss\'e
limit of the class of all finite $L$-structures. The reader may wonder why we did not use this method as well, thus we
will have a look at what we can achieve with Fra\"\i ss\'e limits in our situation. Let us first recall the basic
setup.
\begin{Def}\label{def:fra}
The \emph{age} of a structure~$M$ is the class of all finitely generated structures embeddable in~$M$.

A \emph{Fra\"\i ss\'e class} is a class $\mathcal K$ of finitely generated countable structures satisfying the
following conditions:
\begin{itemize}
\item $\mathcal K$ contains only countably many structures up to isomorphism.
\item Hereditary property (HP): if a finitely generated structure $B$ embeds in $A\in\mathcal K$, then $B\in\mathcal
K$.
\item Joint embedding property (JEP): for any finite set $\{B_0,\dots,B_{n-1}\}\sset\mathcal K$, there exists
$A\in\mathcal K$ such that each $B_i$ embeds in~$A$.
\item Amalgamation property (AP): for any $C,B_0,B_1\in\mathcal K$ and embeddings $g_i\colon C\to B_i$ ($i=0,1$), there
exists $A\in\mathcal K$ and embeddings $f_i\colon B_i\to A$ ($i=0,1$) such that $f_0\circ g_0=f_1\circ g_1$.
\end{itemize}
Note that JEP is equivalent to its special cases $n=0$ and $n=2$. The former amounts to $\mathcal K\ne\nul$; the latter
is, under the assumption of AP and HP, equivalent to the simpler property that $\mathcal K$ contains only one
$0$-generated structure up to isomorphism.

Note also that if the language is countable, finitely generated structures are automatically countable.

A structure~$M$ is \emph{ultrahomogeneous} if for every finitely generated substructures $A,B\sset M$ and every
isomorphism $f\colon A\simeq B$, there exists an automorphism $g$ of~$M$ such that $f\sset g$.
\end{Def}
\begin{Prop}\label{prop:fra}
A class~$\mathcal K$ is a Fra\"\i ss\'e class if and only if it is the age of a countable ultrahomogeneous
structure~$M$. In that case, $M$ is unique up to an isomorphism; it is called the \emph{Fra\"\i ss\'e limit}
of~$\mathcal K$.
\noproof\end{Prop}

Model completions can then be conveniently constructed using Fra\"\i ss\'e limits as follows.
\begin{Prop}\label{prop:fra-mc}
Let $T$ be a universal theory in a finite language consisting of relations and constants. If the class of
finite models of~$T$ has AP and JEP, then it is a Fra\"\i
ss\'e class; its Fra\"\i ss\'e limit~$M$ is a unique countable e.c.\ model of~$T$ up to isomorphism. The theory of~$M$
is the model completion of~$T$.
\noproof\end{Prop}

This works well for $\ec_L$ when $L$ is a finite language with relations and (by considering separately each
quantifier-free diagram) constants. However, it is not applicable if $L$ contains proper functions; this cannot be
circumvented by somehow encoding the structures in a relational language, as the theories have fundamentally different
properties: in particular, model completions constructed by Proposition~\ref{prop:fra-mc} are always
\emph{$\omega$-categorical}, and this seems to be inherent in the method.

Thus, it seems Fra\"\i ss\'e limits are not helpful for showing the existence of~$\ec_L$ in general. Nevertheless, we may still
wonder if limits of suitable Fra\"\i ss\'e classes could provide some interesting models of~$\ec_L$.

Notice that in order to have JEP, whatever class of models we consider must satisfy the same quantifier-free sentences;
in view of quantifier elimination of~$\ec_L$, this amounts to choosing a completion of~$\ec_L$. So, let $L$ be a
countable language, and $T$ a complete extension of~$\ec_L$. We can write $T=\ec_L+\diag(M_0)$, where $M_0$ is a
(possibly empty) structure whose every element is the value of a closed term, i.e., $M_0$ is $0$-generated; we denote $T_0=\diag(M_0)$. There are two
candidate Fra\"\i ss\'e classes of models of~$T_0$ that immediately spring to mind:
\begin{itemize}
\item
The class $\mathcal K_\mathrm{fg}$ of finitely generated models of~$T_0$.
\item
The class $\mathcal K_\mathrm{fin}$ of finite models of~$T_0$.
\end{itemize}
It is easy to verify that $\mathcal K_\mathrm{fg}$ satisfies HP, AP, and JEP. Thus, it is a Fra\"\i
ss\'e class iff it contains countably many structures up to isomorphism, which happens iff $T$ is \emph{small}. (We
already characterized when $T$ is small in Proposition~\ref{prop:categ}.) We leave it to the reader to check that in this case,
the Fra\"\i ss\'e limit of $\mathcal K_\mathrm{fg}$ is the countable saturated model of~$T$.

Clearly, $\mathcal K_\mathrm{fin}$ is nonempty only if $M_0$ is finite, hence we will assume this for the moment. Again,
it is easy to see $\mathcal K_\mathrm{fin}$ has HP, AP, and JEP, hence it is a Fra\"\i ss\'e class iff it is countable
up to isomorphism. This holds iff $L$ contains only finitely many nonnullary symbols; since the finiteness of~$M_0$
implies there are only finitely many constants up to equality in~$T_0$, and nullary relations are in~$T_0$ fixed to true
or false, we can as well assume without loss of generality that $L$ is finite. The Fra\"\i ss\'e limit of $\mathcal
K_\mathrm{fin}$ is then described by the following result.
\begin{Prop}\label{prop:fraisse}
Let $T=\ec_L+\diag(M_0)$, where $M_0$ is a $0$-generated $L$-structure, and $L$ and $M_0$ are \emph{finite}.

Then $T$ has a prime (equivalently: countable atomic) model, which can be characterized as the unique locally
finite (i.e., such that finitely generated submodels are finite) countable model of~$T$, and it can be constructed as the Fra\"\i
ss\'e limit of the class of all finite $L$-structures that extend~$M_0$.
\end{Prop}
\begin{Pf}
By the preceding discussion, $\mathcal K_\mathrm{fin}$ is a Fra\"\i ss\'e class, hence it has a Fra\"\i ss\'e
limit~$M$. Clearly, $M\Sset M_0$ is countable, and locally finite. By general properties of Fra\"\i ss\'e limits (see
e.g.~\cite{hodges:sh}), $M$ is existentially closed in the class of locally finite structures. This in fact implies that $M$ is e.c.\ in the class of all
structures, hence $M\model T$: if $\fii(\vc u)$ is an $\exists_1$ formula, $\vc u\in M$, and $M\sset N\model\fii$, let
\[N'=M\cup\{t^N(\vc u):\text{$t$ is a subterm of $\fii$}\},\]
and make it an $L$-structure by preserving the values of all relations and functions in~$N$ where possible, and
$f^{N'}(\vc v)=a$ for some fixed $a\in M$ if $f^N(\vc v)\notin N'$. Then $M\sset N'\model\fii(\vc u)$, and every
finitely generated submodel of~$N'$ is included in $N'\bez M$ plus a finitely generated submodel of~$M$, hence $N'$ is
locally finite. Thus, $M\model\fii(\vc u)$.

Every locally finite $M'\model T$ is atomic: let $\vc a\in M'$, and $A$ be the submodel of~$M'$ generated
by~$\vc a$. Write $A\bez\{\vc a\}=\{b_i:i<m\}$. Then $\tp(\vc a)$ is generated by $\exists
y_0,\dots,y_{m-1}\,\diag(A)$, where $\diag(A)$ is written using variables $\vc x,\vc y$ in place of $\vc a,\vc
b$.

Thus, $M$ is a countable atomic, hence prime, model of~$T$, and by uniqueness of prime models, every countable
locally finite model of~$T$ is isomorphic to~$M$.
\end{Pf}

There are other cases when $T$ has a prime model: for example, if $T$ is small (see Proposition~\ref{prop:categ}). In fact, we
can give a full description (for countable languages). Recall that a complete countable theory has a prime model if and only
if it is atomic.
\begin{Prop}\label{prop:prime}
Let $T=\ec_L+\diag(M_0)$, where $L$ is countable, and $M_0$ is a $0$-generated $L$-structure. Then $T$ has a prime
model iff it falls in of the following cases:
\begin{enumerate}
\item\label{item:38}
$M_0\model\ec_L$.
\item\label{item:39}
$L$ contains only finitely many nonnullary symbols, all of which are unary.
\item\label{item:40}
$M_0$ is finite, and $L$ contains only finitely many nonnullary symbols.
\end{enumerate}
\end{Prop}
\begin{Pf}
Clearly, \ref{item:38} implies that $M_0$ is a prime model of~$T$. We constructed a prime model in case~\ref{item:40}
in Proposition~\ref{prop:fraisse}.

Assume~\ref{item:39} holds; we will show $T$ is atomic. Let $\fii(\vc x)$ be a $T$-consistent formula, which we may
assume to be quantifier-free. By the argument in the proof of Proposition~\ref{prop:fraisse}, $\fii$ is satisfiable in a model
$M\Sset M_0$ such that $M\bez M_0$ is finite. Let $M_0\sset M\model\fii(\vc a)$ be such that $\lh{M\bez M_0}$ is
minimal possible. Using minimality, we can choose for each $b\in M\bez M_0$ a term $t_b(\vc x)$ (in fact, a subterm
of~$\fii$) such that $t_b(\vc a)=b$. Let $C$ be the (finite) set of all elements of~$M_0$ that are values of subterms
of $\fii(\vc a)$. For each element $b\in M_0$, let us fix a constant term $t_b$ whose value is~$b$. Let $\psi(\vc x)$
be the conjunction of the following formulas:
\begin{itemize}
\item $t_u(\vc x)\ne t_v(\vc x)$ for each $u\in M\bez M_0$, and $u\ne v\in(M\bez M_0)\cup C$;
\item $x_i=t_{a_i}(\vc x)$ for each~$i$;
\item $F(t_u(\vc x))=t_{F(u)}(\vc x)$ for each $u\in M\bez M_0$, and $F\in L$ a nonconstant function symbol;
\item $R(t_u(\vc x))$ or $\neg R(t_u(\vc x))$ (whichever is satisfied by $\vc a$) for each $u\in M\bez M_0$, and
$R\in L$ a nonnullary relation symbol.
\end{itemize}
By construction, $M\model\psi(\vc a)$, hence $\psi$ is consistent. By induction on the length of~$s$, we see that
$T_0\vdash\psi(\vc x)\to s(\vc x)=t_{s(\vc a)}(\vc x)$ for each subterm $s$ of~$\fii$; it follows easily that
$T_0\vdash\psi(\vc x)\to\fii(\vc x)$. We claim that $\psi$ is an atom; by quantifier elimination, it suffices to show
that it implies the quantifier-free type of $\vc a$. Thus, let us consider a model $M'\model\psi(\vc a')$. We may
assume $M'$ is generated by~$\vc a'$. The conjuncts of~$\psi$ ensure that the mapping $b\mapsto t_b(\vc a')$ is a
homomorphism $f\colon M\to M'$ such that $f(a_i)=a'_i$. It is the identity on~$M_0$, and injective on $M\bez M_0$.
Since $\psi$ implies~$\fii$, we have $M'\model\fii(\vc a')$. By the minimality of~$M$, this implies $\lh{M'\bez
M_0}\ge\lh{M\bez M_0}$, hence $f$ must map $M\bez M_0$ to $M'\bez M_0$. Thus, $f$ is in fact an isomorphism of $M$
to~$M'$. 

On the other hand, assume that none of \ref{item:38}--\ref{item:40} holds. Since
$M_0$ does not validate some axiom of~$\ec_L$ as in Definition~\ref{def:ec}, and all elements of~$M_0$ can be denoted by
constant terms, there exists a consistent quantifier-free formula in one variable $\theta(x)$ which is not satisfiable
in~$M_0$. Assuming for contradiction that $T$ is atomic, we may choose $\theta(x)$ to be an atom.

If $L$ includes infinitely many nonnullary relation symbols, let $R(\vc x)$ be one that does not appear in~$\theta$.
Then $\theta(x)\land R(x,\dots,x)$ and $\theta(x)\land\neg R(x,\dots,x)$ are both consistent, contradicting $\theta$
being an atom: taking an arbitrary model $M\model\theta(a)$, we may flip the value of $R(a,\dots,a)$ without
affecting $M\model\theta(a)$. By a similar argument, we obtain a contradiction if $L$ contains infinitely many
nonconstant function symbols.

The remaining case is that $M_0$ is infinite, and $L$ includes an at least binary symbol, say, a relation symbol
$R(\vc x,y)$ (the case of a function symbol is similar). Fix a model $M\model\theta(a)$. Since $M_0$ is infinite, we
can find $c\in M_0$, denoted by a closed term $t$, such that $c$ is not the value of any subterm of $\theta(a)$. Then
we can flip the value of $R(a,\dots,a,c)$ without affecting $M\model\theta(a)$, hence the formulas $\theta(x)\land
R(x,\dots,x,t)$ and $\theta(x)\land\neg R(x,\dots,x,t)$ are both consistent, a contradiction.
\end{Pf}

\section{More on representation}\label{sec:more-representation}

For completeness, let us present a few counterexamples to possible strengthenings of some of the basic claims in
Section~\ref{sec:repr-recurs-funct}.

First, we mentioned that it is enough to represent a specific dprp, or prf, in
a theory in order to show its essential undecidability. In contrast, we will prove that any finite set of trf and rp can be represented in a decidable theory. More generally, it holds even for infinite families
of such functions and predicates as long as they are \emph{uniformly recursive}: here, we call a sequence
$\{F_n:n\in\omega\}$ of recursive functions $F_n\colon\N^{k_n}\to\N$ uniformly recursive if the functions $n\mapsto
k_n$ and $\p{n,w}\mapsto F_n((w)_0,\dots,(w)_{k_n-1})$ are recursive, and similarly for sequences of predicates.
\begin{Prop}\label{prop:unif-rec-dec}
Let $\mathcal R$ be a uniformly recursive sequence of trf and rp. Then there exists a consistent
decidable theory~$T$ and a recursive numeral sequence~$\sigma$ such that $\mathcal R$
is representable in~$T$ w.r.t.~$\sigma$.
\end{Prop}
\begin{Pf}
Let $L$ be the (recursive) language $L_{\mathcal R}$ from Definition~\ref{def:tprf}, and $\p{\N,\mathcal R}$ be the ``standard
$L$-model'' with domain~$\N$ and the $L$-symbols realized by the corresponding elements of~$\mathcal R$. Then $\ec_L+\diag(\p{\N,\mathcal R})=\ec_L+\trep{\mathcal R}$ is a
decidable complete theory by Corollary~\ref{cor:dec}, and it represents $\mathcal R$ w.r.t.\ the sequence $n\mapsto\num n$.
\end{Pf}

The second example serves two-fold purpose. For one, it exhibits that a recursively axiomatizable (or even decidable)
theory may represent a non-recursive predicate w.r.t.\ a non-recursive sequence of numerals. Second, it shows that
representation of rp, or even dprp, does not imply essential undecidability if the
sequence of numerals is not recursive. We will first prove a simple version applying to finite languages.
\begin{Prop}\label{prop:nonrecnumseq}
Let $\mathcal P$ be a finite set of predicates and disjoint pairs. Then there exists a consistent
decidable theory~$T$ that represents $\mathcal P$ (w.r.t.\ a possibly non-recursive sequence of numerals).
\end{Prop}
\begin{Pf}
Since we do not require the elements of~$\mathcal P$ to be recursive to begin with, we may as well extend each disjoint
pair to a predicate, thus we will assume $\mathcal P$ consists of predicates without loss of generality.
Let $L$ be the finite relational language corresponding to~$\mathcal P$. Since $\ec_L$ is decidable, it has a
recursive model~$M$ (using the standard Henkin completion procedure; in fact, in this case, it is not difficult to construct the model explicitly). Let $T=\ec_L+\diag(M)$, which is a decidable complete theory. Since $\ec_L$ is
$\omega$-categorical, $M$ is (a recursive presentation of) its unique countable model, and the countable structure
$\p{\N,\mathcal P}$ embeds in~$M$; let us fix such an embedding $\sigma\colon n\mapsto\num n$, where the elements $\num
n\in M$ are identified with the corresponding constants in the language of~$T$. Then all $\mathcal P$-predicates are
represented in~$T$ w.r.t.~$\sigma$.
\end{Pf}

We could handle countable sets of predicates and disjoint pairs of \emph{bounded arity} with a bit of preprocessing,
but we will need more work to take care of the general case: in particular, note that $\ec_L$ is no longer
$\omega$-categorical (or even small) if $L$ is an infinite relational language; we will use a slightly different theory
instead.

Let $L_0$ be the language that includes an $n$-ary relation $R_n(x_0,\dots,x_{n-1})$ for every $n\ge1$, and $T_0$ be the
universal $L_0$-theory axiomatized by
\[R_n(\vc x)\to\ET_{i<j<n}x_i\ne x_j\]
for each~$n$. The class $\mathcal K$ of finite models of~$T_0$ is easily seen to have HP, AP, and JEP. Crucially, it
contains only countably many nonisomorphic structures: in fact, for any~$n$, there are only finitely many
$T_0$-structures of size~$n$, as the axioms force all but the first $n$ $L_0$-relations to be empty.

Thus, $\mathcal K$ is a Fra\"\i ss\'e class, and it has a Fra\"\i ss\'e limit~$M$. Since any isomorphism between finite
submodels of $M$ extends to an automorphism of~$M$, and there are finitely many isomorphism types of such submodels of
fixed size, it follows that the theory $T^*$ of~$M$ is $\omega$-categorical, and has elimination of quantifiers. Since
all finite(ly generated) models of $T_0$ embed in~$M$, $T^*$ is in fact the model completion of~$T_0$.

It is not difficult to explicitly axiomatize $T^*$ by suitable extension axioms, thus $T^*$ is recursively
axiomatizable, and in fact, decidable (as it is complete). Moreover, $M$ may be presented as a recursive structure.
Thus, its elementary diagram $T^*_M=T^*+\diag(M)$ is a decidable theory.

We will now show that \emph{every} predicate and disjoint pair (with no recursivity assumption) is representable
in~$T^*_M$, using various (mostly nonrecursive) sequences of numerals. In fact, any countable set of predicates and
disjoint pairs can be represented w.r.t.\ the \emph{same} sequence of numerals:
\begin{Prop}\label{prop:nonrec-all}
Every countable set~$\mathcal P$  of predicates and disjoint pairs is representable in the decidable theory~$T^*_M$.
\end{Prop}
\begin{Pf}
As in the proof of Proposition~\ref{prop:nonrecnumseq}, we will assume $\mathcal P$ to consist of predicates. Without loss of
generality, we may also assume that $\mathcal P$ is closed under identification of variables: i.e., if $P\in\mathcal P$
is $n$-ary, and $i<j<n$, the $(n-1)$-ary predicate $P_{i,j}$ defined by
\[P_{i,j}(x_0,\dots,x_{n-2})\iff P(x_0,\dots,x_{j-1},x_i,x_j,\dots,x_{n-2})\]
is in~$\mathcal P$. Let us also define
\[P'(x_0,\dots,x_{n-1})\iff P(x_0,\dots,x_{n-1})\land\ET_{i<j<n}x_i\ne x_j.\]
Then we can reconstruct $P$ from $P'$ and predicates of smaller arity as
\[P(x_0,\dots,x_{n-1})\iff P'(x_0,\dots,x_{n-1})\lor\LOR_{i<j<n}\bigl(x_i=x_j\land
P_{i,j}(x_0,\dots,x_{j-1},x_{j+1},\dots,x_{n-1})\bigr),\]
hence by induction on the arity, we see that all $\mathcal P$-predicates are quantifier-free definable from predicates
from 
$\mathcal P'=\{P':P\in\mathcal P\}$. Thus, without loss of generality, we may replace $\mathcal P$ with~$\mathcal P'$,
i.e., we may assume that all predicates $P\in\mathcal P$ satisfy
\[P(\vc x)\to\ET_{i<j}x_i\ne x_j.\]
We may also assume $\mathcal P$ contains predicates of every arity~$n\in\N$.

For each $n\in\N$, let us fix an enumeration $\{P_{n,m}:m\in\N\}$ of all $n$-ary predicates from~$\mathcal P$ such that
each such predicate occurs more than $n$~times in the enumeration. Finally, for each $n\ge1$, we define
\[P_n(x_0,\dots,x_{n-1})\iff P_{n-1,x_0}(x_1,\dots,x_{n-1})\land\ET_{i>0}x_i\ne x_0.\]

The structure $\p{\N,P_n:n\ge1}$ is a countable model of~$T_0$, hence it embeds in~$M$. Let us fix such an embedding
$\sigma\colon n\mapsto\num n$, where we identify $\num n\in M$ with the corresponding constant in the language
of~$T^*_M$. Thus, $\sigma$ serves as a sequence of numerals, and each relation $P_n$ is represented in~$T^*_M$
w.r.t.~$\sigma$ by the formula $R_n(\vc x)$.

It follows that each $P\in\mathcal P$ is represented in~$T^*_M$ as well: if $P$ is $n$-ary, let us fix
distinct $m_0,\dots,m_n\in\N$ such that $P=P_{n,m_0}=\dots=P_{n,m_n}$. By the pigeonhole principle, all the numbers
$m_0,\dots,m_n$ cannot simultaneously appear in any tuple satisfying $P$; thus, the formula
\[\LOR_{j\le n}R_{n+1}(\num{m_j},x_0,\dots,x_{n-1}).\]
represents $P(x_0,\dots,x_{n-1})$ in~$T^*_M$ w.r.t.~$\sigma$.
\end{Pf}

The main result of this paper shows that representability of prf does not imply
interpretability of~$R$. Another problem in a similar vein is to clarify the relationship between representability of
different types of recursive objects. Specifically, let us consider the following conditions on a theory~$T$:
\begin{enumerate}
\item\label{item:65}
The set $\srp$ is representable in~$T$.
\item\label{item:64}
The set $\sdprp$ is representable in~$T$.
\item\label{item:66}
The set $\srp\cup\{\sucf\}$ is representable in~$T$.
\item\label{item:67}
The set $\sdprp\cup\{\sucf\}$ is representable in~$T$.
\item\label{item:63}
The set $\strf$ is representable in~$T$.
\item\label{item:62}
The set $\sprf$ is representable in~$T$.
\end{enumerate}
We discuss separately the cases with the successor function included because of pathologies exhibited by
representation of predicates and disjoint pairs w.r.t.\ potentially non-recursive sequences of numerals, as seen in
Proposition~\ref{prop:nonrec-all}.

As we already mentioned in Section~\ref{sec:repr-recurs-funct}, it is easy to see that
\ref{item:62}\txto\ref{item:63}, 
\ref{item:67}\txto\ref{item:66}, 
\ref{item:64}\txto\ref{item:65},
\ref{item:62}\txto\ref{item:67}\txto\ref{item:64}, and
\ref{item:63}\txto\ref{item:66}\txto\ref{item:65}.
We wish to show now that no other implications between these six conditions hold in
general. This turns out not to be quite true---once again due to pathologies exhibited by non-recursive sequences of
numerals. In the spirit of battering reality vigorously until it complies with our preformed expectations, we fix this by
considering a more strict notion of implication between representability of classes $\mathcal R_0$ and~$\mathcal R_1$,
namely: does representability of $\mathcal R_0$ in~$T$ w.r.t.\ a sequence of numerals~$\sigma$ imply the
representability of~$\mathcal R_1$ w.r.t.~$\sigma$? This leads to the desired answer.
\begin{Prop}\label{prop:rep-ctrex}
\ \begin{enumerate}
\item\label{item:68}
There exists a theory that represents $\sdprp\cup\{\sucf\}$, but does not represent~$\strf$.
\item\label{item:69}
There exists a theory that represents~$\strf$, but does not represent $\sdprp\cup\{\sucf\}$.
\item\label{item:70}
There exists a theory that represents~$\sdprp$, but does not represent $\sucf$.
\item\label{item:71}
Every theory that represents~$\srp$ also represents arbitrary countable sets of
predicates and disjoint pairs.
\item\label{item:72}
There exists a theory that represents~$\strf$ w.r.t.\ a sequence of
numerals~$\sigma$, but does not represent~$\sdprp$  w.r.t.~$\sigma$.
\end{enumerate}
\end{Prop}
\begin{Pf}
If $\mathcal R\sset\sprf\cup\sdprp$, let $\ec\trep{\mathcal
R}=\ec_{L_{\mathcal R}}+\trep{\mathcal R}$, where $L_{\mathcal R}$ is the language of $\trep{\mathcal R}$.

\ref{item:68}: Let $\mathcal R=\sdprp+\{\sucf\}$. Then the theory
$\ec\trep{\mathcal R}$ represents~$\mathcal R$, and it is supersimple (hence $\prpty{NTP_2}$) by Theorem~\ref{thm:class}.
In contrast, any theory representing trf is $\prpty{TP_2}$ by Proposition~\ref{prop:tp2}.

\ref{item:69}: The (complete) theory $T=\ec\trep\strf$ represents trf. Let $P=\p{P^+,P^-}$ be a recursively inseparable
dprp. We claim that $T$ cannot represent $\{P,\sucf\}$, i.e., it does not interpret $\trep{P,\sucf}$: if it did, then
$\trep{P,\sucf}$ would be also interpretable in a finite-language fragment $T_0$ of~$T$, which would make $T_0$ an
essentially undecidable theory. However, $T_0$ is of the form $\ec\trep{\mathcal R}$ for a finite $\mathcal
R\sset\strf$, thus it is decidable as in the proof of Proposition~\ref{prop:unif-rec-dec}.

\ref{item:70}:
The theory $T=\ec\trep\sdprp$ represents dprp. Assume for contradiction that it represents $\sucf$, i.e., it interprets
$\trep\sucf$. Then $\trep\sucf$ is interpretable in a finite-language fragment $T_0$ of~$T$, which is an extension
of $\ec_L$ for $L$ a finite language with relations and constants. Thus, (any completion of) $T_0$ is
$\omega$-categorical. It follows that $\trep\sucf$ also has an $\omega$-categorical extension, but this is impossible,
as it has infinitely many definable constants.

\ref{item:71}: The argument from the proof of Proposition~\ref{prop:nonrec-all} shows that every
countable set $\mathcal P$ of predicates and disjoint pairs is representable (w.r.t.\ a suitable sequence of numerals)
in any theory that represents all predicates definable in the model~$M$.

\ref{item:72}: The theory $T=\ec\trep\strf$ represents trf w.r.t.\ the sequence of numerals $\num n$. If it also
represented $\sdprp$ w.r.t.\ the same numeral sequence, it would in fact represent $\sdprp\cup\{\sucf\}$, which we
already know to be impossible from~\ref{item:69}.
\end{Pf}

\begin{Rem}\label{rem:shoen}
The theory $T=\ec\trep\strf$ we used in the proof of Proposition~\ref{prop:rep-ctrex} is not recursively axiomatized. Using an
elaborate enumeration of recursive functions, Shoenfield~\cite{shoen:undec-creat} constructed a recursively
axiomatizable theory~$T$ in which all unary trf are representable w.r.t.\ a recursive sequence of numerals~$\sigma$,
but no non-recursive set and no recursively inseparable dprp is representable in~$T$ w.r.t.~$\sigma$.
\end{Rem}

\section*{Acknowledgements}

The question that led to this paper arose from a fruitful discussion with Albert Visser; I am also grateful to him for
pointing me towards Shoenfield~\cite{shoen:undec-creat} (see Remark~\ref{rem:shoen}). I would like to thank Nick Ramsey for
clarification of his work and other comments on the topic. Last but not least, I want to thank the anonymous reviewer
for their useful suggestions and corrections.

The research was supported by grant IAA100190902 of GA AV \v CR, Center of Excellence CE-ITI under the grant
P202/12/G061 of GA \v CR, and RVO: 67985840.

\bibliographystyle{mybib}
\bibliography{mybib}
\end{document}